\documentclass{article}
\pdfoutput=1
\usepackage{graphicx}
\usepackage{latexsym,amsmath,amsfonts,amscd, amsthm}
\usepackage{epsfig}
\usepackage{changebar}
\usepackage{color}
\usepackage{bm}
\usepackage{tikz}

\usepackage{subfigure}
\usepackage{amssymb,amsthm,mathrsfs,dsfont}
\usepackage{changebar}
\usepackage{indentfirst}
\usepackage{verbatim}
\usepackage{amssymb}

\usepackage{multirow}
\usepackage{epstopdf}
\newtheoremstyle{plainNoItalics}{}{}{\normalfont}{}{\bfseries}{.}{ }{}

\PassOptionsToPackage{normalem}{ulem}
\usepackage{ulem}

\numberwithin{equation}{section}

\def\Box{\mbox{ }\rule[0pt]{1.5ex}{1.5ex}}

\begin{document}

\baselineskip=2pc

\theoremstyle{plain}

\newtheorem{algorithm}{ALGORITHM}
\newtheorem{remark}{REMARK}

\newtheorem{theorem}{THEOREM}

\newtheorem{proposition}{PROPOSITION}
\newtheorem{example}{EXAMPLE}


\title{A conservative semi-Lagrangian HWENO method for the Vlasov equation\footnote{Research
was partially supported by NSFC grants 91230110, 11328104, 11571290, 91530107 and NSF DMS-1217008 and DMS-1522777.}}
\author{ Xiaofeng Cai\footnote{School of Mathematical Sciences, Xiamen
University, Xiamen, Fujian, 361005, P.R. China.
E-mail: xfcai89@126.com.} \
 Jianxian Qiu\footnote{School of
Mathematical Sciences and Fujian Provincial Key Laboratory of Mathematical
Modeling \&
High-Performance Scientific Computing, Xiamen University, Xiamen, Fujian,
361005, P.R. China.
E-mail: jxqiu@xmu.edu.cn.}  \
Jingmei Qiu\footnote{Department of Mathematics, University of Houston, Houston, 77204. E-mail: jingqiu@math.uh.edu.}
 }

\maketitle
{\bf Abstract}:
In this paper, we present a high order conservative semi-Lagrangian (SL) Hermite weighted essentially non-oscillatory (HWENO) method for the Vlasov equation based on dimensional splitting [Cheng and Knorr, Journal of Computational Physics, 22(1976)].
The major advantage of HWENO reconstruction, compared with the original WENO reconstruction, is compact.
For the split one-dimensional equation, to ensure local mass conservation, we propose a high order SL HWENO scheme in a conservative flux-difference form, following the work in [J.-M. Qiu and A. Christlieb, Journal of Computational Physics, v229(2010)].
Besides performing dimensional splitting for the original 2D problem, we design a proper splitting for equations of derivatives to ensure local mass conservation of the proposed HWENO scheme.
The proposed fifth order SL HWENO scheme with the Eulerian CFL condition has been tested to work well in capturing filamentation structures without introducing oscillations. We introduce WENO limiters to control oscillations when the time stepping size is larger than the Eulerian CFL restriction.  We perform classical numerical tests on rigid body rotation problem, and demonstrate the performance of our scheme via the Landau damping and two-stream instabilities when solving the Vlasov-Poisson system.


\section{Introduction}
This paper focuses on a high order conservative semi-Lagrangian scheme with high order HWENO reconstruction for the Vlasov-Poisson (VP) simulations based on dimensional splitting. The VP system, arise from collisionless plasma applications, reads as following,
\begin{equation}
\frac{\partial f}{\partial t} + \mathbf{v}\cdot\nabla_{\mathbf{x}}f + \mathbf{E}(t,\mathbf{x})\cdot\nabla_{\mathbf{v}}f=0,
\label{vlasov1}
\end{equation}
and
\begin{equation}
\mathbf{E}(t,\mathbf{x}) = - \nabla_{\mathbf{x}}\phi(t,\mathbf{x}),\ \ -\Delta_{\mathbf{x}}\phi(t,\mathbf{x}) = \rho(t,\mathbf{x}),
\label{poisson}
\end{equation}
where $\mathbf{x}$ and $\mathbf{v}$ are coordinates in phase space $(\mathbf{x},\mathbf{v})\in \mathbb{R}^3\times\mathbb{R}^3$, $\mathbf{E}$ is the electric field, $\phi$ is the self-consistent electrostatic potential and $f(t,\mathbf{x},\mathbf{v})$ is probability distribution function which describes the probability of finding a particle with velocity $\mathbf{v}$ at position $\mathbf{x}$ at time $t$. The probability distribution function  couples to the long range fields via the charge density, $\rho(t,\mathbf{x}) = \int_{\mathbb{R}^3}f(t,\mathbf{x},\mathbf{v})d\mathbf{v}-1$, where we take the limit of uniformly distributed infinitely massive ions in the background. Equations \eqref{vlasov1} and \eqref{poisson} have been nondimensionalized so that all physical constants are one.

Popular methods in fusion simulations include Lagrangian, semi-Lagrangian and Eulerian methods.  Popular lagrangian methods include the particle-in-cell (PIC) \cite{barnes1983implicit,friedman1991multi,jacobs2006high}, Lagrangian particle methods \cite{nature,evstatiev2013variational}; Eulerian methods include weighted essentially non-oscillatory (WENO) coupled with Fourier collocation \cite{zhou2001numerical}, continuous finite element
methods \cite{fem1,fem2}, Runge-Kutta discontinuous Galerkin methods \cite{DG_VP11,de2012discontinuous,heath2012discontinuous,cheng2013study}.
Each method has its own
advantages and limits. For example, Lagrangian methods are well known for its reasonable low computational cost for
high dimensional problems. However, it suffers from statistical noise due to the initial
sampling of macro-particles.
Eulerian methods offer a good alternative to overcome this lack of precision, but they suffer from 'the curse of dimensionality' and the CFL time step restriction.
Compared with the Eulerian approach,
the SL methods is relief from the CFL time step restriction, because information is being propagated along characteristics.

Among SL schemes, the scheme based on dimensional splitting, introduced by Cheng and Knorr originally \cite{cheng1976integration}, is very popular
SL scheme with high order cubic spline interpolation was proposed in \cite{sonnendrucker1999semi}. A positivity preserving and flux conservative finite volume SL scheme with ENO reconstruction for the VP system is proposed in \cite{filbet2001conservative} and the scheme for the guiding center
Vlasov model is proposed in \cite{crouseilles2010conservative}.
A conservative
finite different semi-Lagrangian scheme with WENO reconstruction is proposed in \cite{slweno}; later the algorithm is generalized to
variable coefficient case \cite{qiu2010conservative,qiuconservative_cicp} and maximum principle preserving limiter  \cite{xiong2014high_vlasov}.
In the finite element discontinuous Galerkin framework, there are SL discontinuous Galerkin schemes with positivty preserving limiters \cite{qiu2011positivity,rossmanith2011positivity} and hybrid SL finite element-finite difference methods in \cite{guo2013hybrid}.
High order propagation methods based on Hermite interpolation are proposed in \cite{nakamura1999cubic,filbet2003comparison,besse2003semi,besse2008convergence,yang2014conservative}.
HWENO scheme was introduced in \cite{qiu2004hermite} and further developed in \cite{zhu2008class,liu2015finite} for hyperbolic conservation laws. Besides the original equation, one also evolves equations of derivatives in the WENO fashion. Hence their reconstruction stencils are more compact than the original WENO scheme \cite{jiang1996efficient}, given the same order of approximation.
A similar technique, called the CIP (Constrained Interpolation Profile/Cubic Interpolated Propagation) scheme \cite{yabe2001constrained}, has also been proposed for the VP system
\cite{nakamura1999cubic}.
 In \cite{besse2003semi}, besides the function values themselves, the gradients of the function are evolved in a semi-Lagrangian fashion. In \cite{yang2014conservative}, a semi-Lagrangian HWENO is proposed without evolving the gradients. The proposed method in this paper successfully couples the semi-Lagrangian framework with HWENO method. Compared with those earlier work, our method achieves local mass conservation, has relatively compact stencil and is able to capture under-resolve solution structures without numerical artifacts.

 In this paper, we design a conservative SL HWENO scheme based on dimensional splitting. Its stencil is more compact than the regular WENO scheme. We follow the idea in \cite{slweno} to express the SL update in a flux difference conservative form. There are several new ingredients we developed in this paper in order to ensure the effectiveness and robustness of the proposed scheme: firstly, we design a SL HWENO scheme in a flux difference form for 1D problem; secondly, we apply  a special splitting similar to the splitting in the CIP method \cite{nakamura1999cubic} for mass conservation for high dimensional problems; thirdly, we introduce WENO limiters for controlling oscillations.

The paper is organized as follows. In Section \ref{slhweno1d}, we introduce a conservative form of high order SL HWENO method for 1D transport problems. In Section \ref{slhweno2d}, we
introduce a conservative SL HWENO for VP system by a special form of splitting, and introduce WENO limiters. In Section \ref{tests}, we
present our numerical results for basic test problems, such as linear advection and rigid body
rotation, and for  the VP simulations. Concluding remarks
are given in Section \ref{conclusion}.

\section{Conservative SL HWENO method for 1D transport problem}
\label{slhweno1d}
In this section, we introduce the SL Hemite interpolation in a flux-difference form for 1D transport problems. Then we incorporate the HWENO mechanism into the flux function reconstruction procedure to realize a non-oscillatory capturing of shocks.
\subsection{The SL Hemite interpolation in a flux-difference form for 1D transport problem}

In this section, we consider a 1D transport problem
\begin{equation}
\label{linear1}
f_t+ vf_x=0,\quad f(x,t=0)=f_0(x),\quad \mbox{on} \quad [a,b],
\end{equation}
where $v$ is a constant.
For the Hermite method, we also consider the evolution equation for the solution's derivative $g \doteq f_x$. For the linear transport problem \eqref{linear1}, $g$ satisfies the same linear transport equation
\[
g_t + v g_x = 0.
\]
We discretize the domain $[a, b]$ as
\begin{equation*}
a = x_{\frac12} < x_{\frac32} < \cdots < x_{N+\frac12} = b,
\end{equation*}
with the uniform grid points $x_i=a+(i-\frac12)\Delta x$ and the cell size $\Delta x = x_{i+\frac12}-x_{i-\frac12}$.
We let $I_i = [x_{i-\frac12},x_{i+\frac12}]$ and $I_{i+\frac12}=[x_i,x_{i+1}]$.
 In the HWENO approach, the numerical solutions associated with each grid point are point values $f^n_i$ and derivatives $g^n_i$. Here the subscript $i$ means the solution at the grid point $x_i$ and the superscript $n$ means the solution at time level $t^n$. To design a SL HWENO scheme, we update $\{f^{n+1}_i, g^{n+1}_i\}_{i=1}^N$ from the corresponding solutions at time $t^n$.

For the linear problem \eqref{linear1} with constant characteristic speed, the solutions $f_i^{n+1}$ and $g_i^{n+1}$ can be obtain by shifting the information at $t^n$ in the SL framework. We define the amount of shift $\frac{v\Delta t}{\Delta x}$ as $xshift$. There are three cases of $xshift$: shift to the right by some amount less than half a cell ($xshift \in [0,\frac12]$), shift to the left by some amount less than half a cell ($xshift\in[-\frac12,0]$) and shift a distance greater than half a cell ($|xshift|> \frac12$).

To illustrate the idea, we only consider Hermite interpolations for $xshift\in[0,\frac12]$, while the one for $xshift\in[-\frac12,0]$ is mirror symmetric with respect to $x_i$ of the previous interpolations.
In the case when $|xshift|>\frac12$, whole grid shifting is carried out and followed by a final update based on the procedure for $xshift\in [-\frac12,\frac12]$. In the following, we present the Hermite interpolation with cubic polynomials. Higher order schemes will be discussed later.
\begin{enumerate}
  \item The underlying function at $t^n$ can be approximated by a Hermite-type reconstruction, based on   the stencil $\{ f_{i-1}^n,f_i^n,g_{i-1}^n,g_i^n \}$,
\begin{eqnarray*}
\widetilde{f}_{i-\frac12}^n(\xi) &= & f_{i}^n -g_{i}^n \Delta x\xi+ \left( 2g_{i}^n \Delta x -3 f_{i}^n + 3f_{i-1}^n +g_{i-1}^n \Delta x \right)\xi^2 \\
  & &+ \left( 2 f_{i}^n - g_{i}^n \Delta x -2 f_{i-1}^n-g_{i-1}^n \Delta x \right)\xi^3, \notag
\end{eqnarray*}
where $\xi(x)=\frac{x-x_i}{x_{i-1}-x_i}\in [0,1],x\in I_{i-\frac12}$.

  \item $f^{n+1}_i$ and $g^{n+1}_i$ can be obtained by tracing the characteristic back to time $t=t^n$ and evaluating the interpolant $\widetilde{f}_{i-\frac12}^n(\xi)$ at the foot of characteristics $x^\star = x_i - v \Delta t$,
\begin{eqnarray}
 f_i^{n+1} &= &f_i^n - \xi_0 ( (3 f_i^n \xi_0 -2f_i^n\xi_0^2)    -  (  3f_{i-1}^n \xi_0 -2f_{i-1}^n \xi_0^2  )) \notag \\
       &&-g_i^n \Delta x\xi_0 +(2g_i^n \Delta x+g_{i-1}^n \Delta x)\xi_0^2 +(-g_i^n \Delta x-g_{i-1}^n\Delta x)\xi_0^3 \label{hermite11}\\
g_i^{n+1} & =& g_i^n + \left(-4g_i^n + \frac{6f_i^n-6f_{i-1}^n }{\Delta x} - 2g_{i-1}^n\right)\xi_0    \notag     \\
      &  &+  \left(-\frac{6f_i^n-f_{i-1}^n }{\Delta x} + 3g_i^n + 3g_{i-1}^n \right)\xi_0^2,\label{hermite22}
\end{eqnarray}
where $\xi_0 = \frac{x^\star - x_i}{x_{i-1} - x_i}$.
\end{enumerate}

For the above linear scheme \eqref{hermite11} and \eqref{hermite22}, we have the following mass conservation result.
\begin{proposition}
\label{proposition1}
If $\mathop{\sum_{i=1}^N} g_i^n\equiv 0$ and with periodic boundary condition, then the scheme \eqref{hermite11} and \eqref{hermite22} conserve the total mass, i.e., $\mathop{\sum_{i=1}^N} f_i^{n+1}\equiv\mathop{\sum_{i=1}^N} f_i^{n}$ and $\mathop{\sum_{i=1}^N} g_i^{n+1}\equiv 0$.
\end{proposition}
{\bf Proof:}
\begin{eqnarray*}
\sum_{i=1}^N f_i^{n+1} &=& \sum_{i=1}^N [f_i^n - \xi_0 ( (3 f_i^n \xi_0 -2f_i^n\xi_0^2)    -  (  3f_{i-1}^n \xi_0 -2f_{i-1}^n \xi_0^2   ))\\
       & &-g_i^n \Delta x\xi_0 +(2g_i^n \Delta x+g_{i-1}^n \Delta x)\xi_0^2 +(-g_i^n \Delta x-g_{i-1}^n\Delta x)\xi_0^3]   \\
       & =&\sum_{i=1}^N f_i^n -\xi_0(  (3 f_N^{n} \xi_0 -2f_N^n\xi_0^2) -(3 f_0^n \xi_0 -2f_0^n\xi_0^2) )  \\
       & &+ \sum_{i=1}^N g_{i-1}^n\Delta x \xi_0^2-  \sum_{i=1}^Ng_{i-1}^n\Delta x\xi_0^3  \\
       & = &\sum_{i=1}^N f_i^n  \ \text{(periodic boundary condition)},
\end{eqnarray*}
       \begin{eqnarray*}
\sum_{i=1}^N g_i^{n+1} & = & \sum_{i=1}^N\left[g_i^n + \left(-4g_i^n + \frac{6f_i^n-6f_{i-1}^n }{\Delta x} - 2g_{i-1}^n\right)\xi_0          \right. \\
      & & \left.  +  \left(-\frac{6f_i^n-f_{i-1}^n }{\Delta x} + 3g_i^n + 3g_{i-1}^n \right)\xi_0^2\right]\\
      & = &\sum_{i=1}^Ng_i^n+\frac{6f_N^n-6f_{0}^n }{\Delta x}\xi_0-\frac{6f_N^n-f_{0}^n }{\Delta x}\xi_0^2 +2\sum_{i=1}^N g_{i-1}^n\Delta x \xi_0^2-  3\sum_{i=1}^Ng_{i-1}^n\Delta x\xi_0^3  \\
      & = &\sum_{i=1}^Ng_i^n\ \text{(periodic boundary condition)}.
\end{eqnarray*}
Hence the proposition is proved.      \hspace{\stretch{1}} $\Box$

In order to guarantee $\sum_{i=1}^N g_i^0\equiv 0$ in the assumption of the proposition, we introduce a sliding average function $h(x)$ in \cite{shu1989efficient} which satisfies
\begin{equation*}
f(x)=\frac{1}{\Delta x}\int_{x-\frac{\Delta x}{2}}^{x+\frac{\Delta x}{2}} h(\xi) d\xi,
\end{equation*}
then
\begin{equation*}
\label{JMQiu}
g(x)=f(x)_x= \frac{1}{\Delta x }\left( h\left(x+\frac{\Delta x}{2}\right) - h\left( x-\frac{\Delta x}{2} \right)  \right).
\end{equation*}
Thus $\sum_{i=1}^N g_i^0=\sum_{i=1}^N(h_{i+\frac12}^0 - h_{i-\frac12}^0 )\equiv 0$ where $h^0_{i\pm\frac12} \approx h(x
\pm \frac{\Delta x}{2})$ can be obtained by reconstruction from $\{f^0_j\}_j$.

In the following, we will adopt a matrix notation for presentation of the Hermite interpolation. The matrix $A$ will denote the interpolation matrix. We use $A(i,j)$ to denote the element at the $i^{th}$ row and $j^{th}$ column, $A(i,:)$ to denote the $i^{th}$ row of $A$, and $A(:,j)$ to denote the $j^{th}$ column of $A$.

We rewrite the scheme \eqref{hermite11} and \eqref{hermite22} into a flux difference form, in order to ensure local mass conservation, especially when the nonlinear HWENO mechanism is applied. In order to do so, we propose to update $\{f_i^n,h_{i+\frac12}^n\}_i$ instead of $\{f_i^n,g_i^n\}_i$, observing that $g_i^n$ can be recovered from $\{h_{i+\frac12}^n\}_i$ by $g_i^n=\left(h_{i+\frac12}^n-h_{i-\frac12}^n\right)/\Delta x$. Specifically, \eqref{hermite11} can be rewritten in the following flux difference form using the new $\{h_{i+\frac12}^n\}_i$,
\begin{eqnarray}
f_i^{n+1} &=& f_i^n - \xi_0 ( (3 f_i^n \xi_0 -2f_i^n\xi_0^2)    -  (  3f_{i-1}^n \xi_0 -2f_{i-1}^n \xi_0^2   )) \notag\\
       & &-g_i^n \Delta x\xi_0 +(2g_i^n \Delta x+g_{i-1}^n \Delta x)\xi_0^2 +(-g_i^n \Delta x-g_{i-1}^n\Delta x)\xi_0^3  \notag \\
       &= &f_i^n - \xi_0 ( f_i^n(3 \xi_0 -2\xi_0^2)    -  f_{i-1}^n(3 \xi_0 -2\xi_0^2)  ) + g_i^n\Delta x(-\xi_0 +2\xi_0^2 -\xi_0^3) +g_{i-1}^n \Delta x( \xi_0^2-\xi_0^3)  \notag \\
       & =&f_i^n - \xi_0 ( f_i^n(3 \xi_0 -2\xi_0^2)    -  f_{i-1}^n(3 \xi_0 -2\xi_0^2) ) \notag \\
       &  &- \left(h_{i+\frac12}^{n} -h_{i-\frac12}^{n}\right)\xi_0(1 -2\xi_0 +\xi_0^2) -\left(h_{i-\frac12}^{n} -h_{i-\frac32}^{n}\right)\xi_0(- \xi_0+\xi_0^2)  \notag \\
       & =&f_i^n- \xi_0 \left\{ \left[ f_i^n(3\xi_0-\xi_0^2) + h_{i+\frac12}^n(1-2\xi_0+\xi_0^2) +h_{i-\frac12}^n(-\xi_0-\xi_0^2) \right]   \right. \notag\\
      & &\left.- \left[  f_{i-1}^n(3\xi_0-\xi_0^2) + h_{i-\frac12}^n(1-2\xi_0+\xi_0^2) +h_{i-\frac32}^n(-\xi_0-\xi_0^2) \right]  \right\}  \notag \\
      &=&f_i -\xi_0(\widehat{f}_{i+\frac12}^n(\xi_0) -  \widehat{f}_{i-\frac12}^n (\xi_0) ). \notag
       \end{eqnarray}
where
\begin{align*}
\widehat{f}_{i-\frac12}^n (\xi_0) &=   (f_{i-1}^n,h_{i-\frac12}^n ,h_{i-\frac32}^n )\cdot C_3^L \cdot (1,\xi_0,\xi_0^2)^\prime
\end{align*}
with
\begin{equation*}
C_3^L=\left(
\begin{array}{ccc}
0 & 3 & -2\\
1 & -2 & 1 \\
0 & -1 & 1
\end{array}\right).
\end{equation*}
We update $g^{n+1}_i$ by
\begin{equation}
g_i^{n+1} =\frac{ h_{i+\frac12}^{n+1} -h_{i-\frac12}^{n+1} }{\Delta x}
\end{equation}
where
\begin{equation}
\label{hi}
h_{i-\frac12}^{n+1}=(f_{i-1}^n,h_{i-\frac12}^n ,h_{i-\frac32}^n )\cdot D_3^L \cdot (1,\xi_0,\xi_0^2)^\prime
\end{equation}
with
\begin{equation*}
D_3^L=\left(
\begin{array}{ccc}
0 & 6 & -6\\
1 & -4 & 3 \\
0 & -2 & 3
\end{array}\right).
\end{equation*}


\begin{remark}
The flux-difference form for the SL finite difference scheme was originally proposed in \cite{slweno}. There are two main advantages to work with the flux difference form:
  \begin{enumerate}
\item The flux difference form can ensure local mass conservation.
\item We can design a nonlinear HWENO mechanism for the flux reconstructions, see discussions in the next subsection.
\end{enumerate}
In order to work with the flux difference form, it is crucial to work with the $\{h_{i+\frac12}^n\}_i$ instead of the original derivative function $g = f_x$.
\end{remark}

\begin{remark}
We observe that $D_3^L(:,k)=kC_3^L(:,k),\ k=1,2,3$.
\end{remark}

\begin{remark}
The case presented here is for the third order scheme. Similar procedure can be used to obtained higher order scheme, e.g. the fifth order case with HWENO is presented in the next subsection.
\end{remark}

\subsection{HWENO reconstruction for flux functions}
\label{SLHWENO}

In general, high order fixed stencil reconstruction of numerical fluxes performs well when the solution is smooth. However, around discontinuities, oscillations will be introduced. In this subsection, a nonlinear SL HWENO procedure is introduced for reconstructing the flux $\widehat{f}_{i-\frac12}^n(\xi)$. By adaptively assigning nonlinear weights to neighboring candidate stencils, the nonlinear HWENO reconstruction preserves high order accuracy of the linear scheme around smooth regions of the solution, while producing a sharp and essentially non-oscillatory capture of discontinuities.

We adopt the idea of the HWENO reconstruction \cite{qiu2004hermite,liu2015finite} into the proposed conservative SL framework. We present a fifth order HWENO reconstruction as an example. Similar procedure can be generalized to higher order case.

 Our discussion will be focused on the case of $xshift\in[-\frac12,\frac12]$. As before, the case of $|xshift| >\frac12$ will be handled with a whole grid shift followed by the case of $xshift\in[-\frac12,\frac12]$ to account for the fractional remainder.

When $xshift\in[0,\frac12]$, the fifth order conservative SL method the $\{f_{i-2}^n,f_{i-1}^n,f_i^n,f_{i+1}^n,g_{i-2}^n,g_{i+1}^n\}$ is the following,
\begin{eqnarray}
f_i^{n+1} &=& f_i^n + \left( -\frac{8}{27}f_{i-2}^n + f_{i-1}^n - \frac{19}{27}f_{i+1}^n + \frac{2}{9} g_{i+1}^n\Delta x - \frac{1}{9}g_{i-2}^n\Delta x \right)\xi_0 \notag\\
& & + \left(-\frac{1}{18} g_{i-2}^n\Delta x -\frac{2}{9}g_{i+1}^n\Delta x -\frac{7}{4}f_i^n - \frac{19}{108}f_{i-2}^n + f_{i-1}^n + \frac{25}{27}f_{i+1}^n \right) \xi_0^2 \notag\\
& & + \left( \frac{1}{6}g_{i-2}^n\Delta x -\frac{1}{6}g_{i+1}^n\Delta x + \frac14 f_i^n + \frac{5}{12}f_{i-2}^n-\frac{3}{4}f_{i-1}^n + \frac{1}{12}f_{i+1}^n \right) \xi_0^3 \notag\\
& & + \left( \frac{1}{18}g_{i-2}^n\Delta x + \frac{2}{9} g_{i+1}^n\Delta x  + \frac{3}{4}f_i^n + \frac{19}{108}f_{i-2}^n - \frac12f_{i-1}^n - \frac{23}{54}f_{i+1}^n \right) \xi_0^4 \notag\\
& & + \left( -\frac{1}{18}g_{i-2}^n\Delta x -\frac{1}{18}g_{i+1}^n\Delta x - \frac14 f_i^n - \frac{13}{108}f_{i-2}^n + \frac14 f_{i-1}^n + \frac{13}{108} f_{i+1}^n \right) \xi_0^5, \notag
\end{eqnarray}
Using the flux difference form for $g$ function, $g_i^{n} =\frac{ h_{i+\frac12}^{n} -h_{i-\frac12}^{n} }{\Delta x}$ in $t^n$, then
\begin{equation*}
\begin{split}
f_i^{n+1} & =  f_i^n - \xi_0 ( (f_{i-2}^n,f_{i-1}^n,f_i,f_{i+1}^n,h_{i-\frac52}^n,h_{i-\frac32}^n,h_{i+\frac12}^n,h_{i+\frac32}^n) \cdot B_5^L \cdot (1,\xi_0,\xi_0^2,\xi_0^3,\xi_0^4)' )
\end{split}
\end{equation*}
where
\begin{equation*}
B_5^L =
\left(
\begin{array}{ccccc}
\frac{8}{27}  &   \frac{19}{108}   &    -\frac{5}{12}    &    -\frac{19}{108}   &    \frac{13}{108}  \\
   -1         &        -1          &        \frac34      &      \frac12         &      - \frac14     \\
   0          &      \frac74       &      -\frac14       &    -\frac34          &    \frac14         \\
\frac{19}{27} &   -\frac{25}{27}   &    -\frac{1}{12}    &     \frac{23}{54}    &    -\frac{13}{108} \\
  -\frac19    &    - \frac{1}{18}  &      \frac{1}{6}    &      \frac{1}{18}    &    -\frac{1}{18}   \\
  \frac19     &    \frac{1}{18}    &    -\frac{1}{6}     &     -\frac{1}{18}    &     \frac{1}{18}   \\
  \frac29     &   -\frac29         &    -\frac16         &     \frac29          &    -\frac{1}{18}   \\
  -\frac29    &    \frac29         &     \frac16         &     -\frac29         &     \frac{1}{18}
\end{array}
\right).
\end{equation*}

Then $f_i^{n+1}$ can be written in the flux difference form,
\begin{equation}
f_i^{n+1} = f_i^n - \xi_0 ( (f_{i-1}^n,f_i^n,f_{i+1}^n,h_{i-\frac32}^n,h_{i+\frac32}^n)\cdot C_5^L - (f_{i-2}^n,f_{i-1}^n,f_i^n,h_{i-\frac52}^n,h_{i+\frac12}^n)\cdot C_5^L )\cdot (1,\xi_0,\xi_0^2,\xi_0^3,\xi_0^4)'
\label{fflux}
\end{equation}
where
\begin{equation*}
C_5^L =
\left(
\begin{array}{ccccc}
   -\frac{8}{27}   &    -\frac{19}{108}     &    \frac{5}{12}     &    \frac{19}{108}     &   -\frac{13}{108}  \\
    \frac{19}{27}  &     \frac{89}{108}     &   -\frac{1}{3}      &    -\frac{35}{108}    &    \frac{7}{54}    \\
    \frac{19}{27}  &    -\frac{25}{27}      &   -\frac{1}{12}     &     \frac{23}{54}     &   -\frac{13}{108}  \\
    \frac{1}{9}    &     \frac{1}{18}       &   -\frac{1}{6}      &    -\frac{1}{18}      &    \frac{1}{18}    \\
   -\frac{2}{9}    &     \frac{2}{9}        &    \frac{1}{6}      &    -\frac{2}{9}       &    \frac{1}{18}
\end{array}
\right)
\end{equation*}

And we have the flux difference form the derivative $g_i^{n+1}$ in $t^{n+1}$,
\begin{equation*}
\begin{split}
g_i^{n+1}  =  \frac{h_{i+\frac12}^{n+1}  -h_{i-\frac12}^{n+1} }{\Delta x}
\end{split}
\end{equation*}
where
\begin{equation}
h_{i-\frac12}^{n+1} = (f_{i-2}^n,f_{i-1}^n,f_i^n,h_{i-\frac52}^n,h_{i+\frac12}^n)\cdot D_5^L \cdot(1,\xi_0,\xi_0^2,\xi_0^3,\xi_0^4)'
\label{hflux}
\end{equation}
where $D_5^L$ satisfies $D_5^L(:,k)=kC_5^L(:,k),\ k=1,\cdots,5$.

When $xshift\in[-\frac12,0]$, we update $\{f_i^{n+1},h_{i+\frac12}^{n+1}\}_i$ by the following formulas,
\begin{equation}
f_i^{n+1} = f_i^n + \xi_0 ( \widehat{f}_{i+\frac12}^n(\xi_0) - \widehat{f}_{i-\frac12}^n(\xi_0) ),
\end{equation}
where the flux function
\begin{equation}
\widehat{f}_{i-\frac12}^n(\xi) = ( f_{i-1}^n ,f_i^n ,f_{i+1}^n , h_{i-\frac32} , h_{i+\frac32} ) \cdot C_5^R \cdot (1,\xi_0,\xi_0^2 ,\xi_0^3,\xi_0^4)',
\end{equation}
where
\begin{equation}
C_5^R =
\left(
\begin{array}{ccccc}
 \frac{19}{27}  &  -\frac{25}{27}   &   -\frac{1}{12}  &  \frac{23}{54}   &   -\frac{13}{108}  \\
 \frac{19}{27}  &  \frac{89}{108}   &   -\frac13       & -\frac{35}{108}  &    \frac{7}{54}    \\
 -\frac{8}{27}  &  -\frac{19}{108}  &    \frac{5}{12}  &  \frac{19}{108}  &   -\frac{13}{108}  \\
 -\frac29       &   \frac29         &    \frac16       &  -\frac29        &    \frac{1}{18}    \\
 \frac19        &   \frac{1}{18}    &   -\frac16       &  -\frac{1}{18}   &    \frac{1}{18}
\end{array}
\right).
\end{equation}

And
\begin{equation}
h_{i-\frac12}^{n+1}
=
(f_{i-1}^n , f_i^n,f_{i+1}^n,h_{i-\frac32}^n,h_{i+\frac32}^n)\cdot D_5^R \cdot (1,\xi_0,\xi_0^2,\xi_0^3,\xi_0^4)'
\end{equation}
where $D_5^R$ satisfies $D_5^R(:,k)=kC_5^R(:,k),\ k=1,\cdots,5$.

In the following, we illustrate the corresponding HWENO reconstruction of flux functions.
We only discuss the HWENO reconstruction for the flux $\widehat{f}_{i-\frac12}^n$ and $h_{i-\frac12}^{n+1}$ when $xshift\in[0,\frac12]$.
When $xshift\in[-\frac12,0]$, the flux $\widehat{f}_{i-\frac12}^n$ and $h_{i-\frac12}^{n+1}$ could be reconstructed symmetrically with respect to $x_i$. From equations \eqref{fflux} and \eqref{hflux}, the stencil $\{f_{i-2}^n,f_{i-1}^n,f_i^n,h_{i-\frac52}^n,h_{i+\frac12}^n\}$ is used to construct the flux $\widehat{f}_{i-\frac12}^n(\xi)$ and $h_{i-\frac12}^{n+1}$.
It is composed of the information from three potential stencils
\begin{equation}
S_1 = \{h_{i-\frac52}^n,f_{i-2}^n,f_{i-1}^n\},\
S_2 = \{ f_{i-2}^n,f_{i-1}^n,f_i^n \},\
S_3 = \{ f_{i-1}^n,f_i^n,h_{i+\frac12}^n\}.
\end{equation}
Intuitively, in regions where the function is smooth, we want to use information from $S_1 , S_2$ and $S_3$ in an optimal way, to obtain a fifth order approximation.
On  the other hand, around a big jump, we only want to use the information from the relatively smooth stencil.
Following \cite{slweno}, we only use the HWENO mechanism in adaptively reconstructing the coefficients in front of the constant 1 in the equation for $\widehat{f}_{i-\frac12}^n$ and $h_{i-\frac12}^{n+1}$, while leaving coefficients for $\xi_0,\xi_0^2,\xi_0^3,\xi_0^4$ unchanged.
We can observe that the first column of matrix $C_5^L$ is the same as that of $D_5^L$. Thus we only consider the HWENO procedure for constructing $\widehat{f}_{i-\frac12}^n$,
\begin{enumerate}
  \item Compute the linear weights, $\gamma_1,\gamma_2$ and $\gamma_3$, such that
  \begin{equation*}
  \begin{split}
  &(f_{i-2}^n,f_{i-1}^n,f_i^n,h_{i-\frac52}^n,h_{i+\frac12}^n)\cdot C_5^L(:,1)  \\
   &=
  \gamma_1 (f_{i-2}^n,f_{i-1}^n,h_{i-\frac52}^n )\cdot(-2,2,1)'+\gamma_2 (f_{i-2}^n,f_{i-1}^n,f_i^n)\cdot(-\frac16,\frac56,\frac13)' \\
  &+ \gamma_3(f_{i-1}^n,f_i^n,h_{i+\frac12}^n)\cdot(\frac14,\frac54,-\frac12)',
  \end{split}
  \end{equation*}
   where $(f_{i-2}^n,f_{i-1}^n,h_{i-\frac52}^n )\cdot(-2,2,1)',\ (f_{i-2}^n,f_{i-1}^n,f_i^n)\cdot(-\frac16,\frac56,\frac13)'$ and
   $(f_{i-1}^n,f_i^n,h_{i+\frac12}^n)\cdot(\frac14,\frac54,-\frac12)'$ are third order reconstructions of fluxes from three stencils $S_1,S_2$ and $S_3$ respectively. From equation \eqref{fflux}, $\gamma_1=\frac19$, $\gamma_2 =\frac49$ and $\gamma_3 = \frac49$.

  \item We compute the smoothness indicator, denoted by $\beta_j$, for each stencil $S_j$, which measures how smooth the function $p_j(x)$ is in the target cell $I_i$. The smaller this smoothness indicator $\beta_j$, the smoother the function $p_j(x)$ is in the target cell. We use the same recipe for the smoothness indicator as in \cite{jiang1996efficient},
      \begin{equation*}
      \beta_j = \mathop{ \sum_{l=1}^2}
      \int_{I_i}\Delta x^{2l-1}\left( \frac{\partial}{\partial x^l}p_j(x) \right)^2dx.
      \end{equation*}
      In the actual numerical implementation the smoothness indicators $\beta_j$ are written out explicitly as quadratic forms of the points $\{f_i^n,h_{i+\frac12}^n\}_i$ in the stencil,
      \begin{equation*}
      \begin{split}
      \beta_1 &= \frac{13}{3}\left(-\frac94 f_{i-2}^n + \frac32 h_{i-\frac52}^n + \frac34 f_{i-1}^n\right)^2
       +   \left(\frac{31}{4}f_{i-2}^n-\frac92 h_{i-\frac52}^n -\frac{13}{4}f_{i-1}^n\right)^2 ,\\
       \beta_2 & = \frac{13}{12}\left(-f_{i-2}^n +2 f_{i-1}^n -f_i^n \right)^2
       + \left(  -\frac32f_i^n + 2f_{i-1}^n - \frac12f_{i-2}^n  \right)  ^2                                          , \\
       \beta_3 & =\frac{13}{3}\left(-\frac94 f_i^n + \frac34 f_{i-1}^n +\frac32 h_{i+\frac12}^n\right)^2
       +\left(\frac54 f_i^n +\frac14 f_{i-1}^n -\frac32 h_{i+\frac12}^n\right)^2.
      \end{split}
      \end{equation*}

  \item We compute the nonlinear weights based on the smoothness indicators.
  \begin{equation*}
  \omega_j = \frac{\overline{\omega}_j}{  \mathop{\sum_{k=1}^3}\overline{\omega}_k },\ j= 1,2,3,\
  \overline{\omega}_k = \frac{\gamma_k}{\epsilon + \beta_k}
  \end{equation*}
  where $\epsilon$ is a small number to prevent the denominator to becoming zero. In our numerical tests we take $\epsilon$ to be $10^{-6}$.

  \item Compute numerical fluxes constructed in HWENO fashion. Define the matrix $\widetilde{C}_5^L$ and $\widetilde{D}_5^L$ as,
  \begin{equation*}
  \widetilde{D}_5^L(:,1)=\widetilde{C}_5^L(:,1)
  = \omega_1 \cdot(-2,2,0,1,0) + \omega_2 \cdot(-\frac16,\frac56,\frac13,0,0) +\omega_3 \cdot(0,\frac14,\frac54,0,-\frac12)
  \end{equation*}
  \begin{equation*}
  \widetilde{C}_5^L(:,k) = C_5^L(:,k),\
  \widetilde{D}_5^L(:,2) = kC_5^L(:,k),\ k=2,\cdots,5.
  \end{equation*}

  The updated numerical flux is computed using $\widetilde{C}_5^L$ and $\widetilde{D}_5^L$, i.e.,
  \begin{equation}
  \widehat{f}_{i-\frac12}^n (\xi_0)= (f_{i-2}^n,f_{i-1}^n,f_i^n,h_{i-\frac52}^n,h_{i+\frac12}^n)\cdot \widetilde{C}_5^L \cdot (1,\xi_0,\xi_0^2,\xi_0^3,\xi_0^4)',
  \end{equation}
  \begin{equation}
  h_{i-\frac12}^{n+1} = (f_{i-2}^n,f_{i-1}^n,f_i^n,h_{i-\frac52}^n,h_{i+\frac12}^n)\cdot \widetilde{D}_5^L \cdot (1,\xi_0,\xi_0^2,\xi_0^3,\xi_0^4)'.
  \end{equation}
\end{enumerate}

\section{Strang splitting SL HWENO scheme for the VP system}
\label{slhweno2d}

In this section, we extend the SL HWENO scheme in the previous section for solving the 1D VP system.

Denoting by $f(t,x,v)\geq 0$ the distribution function of electrons in phase space and by $E(t,x)$ the self-consistent electric field, the dimensionless VP systems reads as
\begin{equation}
\frac{\partial f}{\partial t} + v \frac{\partial f}{\partial x} + E(t,x)\frac{\partial f}{\partial v}=0,
\label{vlasov}
\end{equation}
\begin{equation}
\frac{dE}{dx}(t,x) = \rho (t,x) = \int_{-\infty}^{+\infty} f(t,x,v) dv -1,
\label{eletric}
\end{equation}
on the domain $[a, b] \times [-L, L]$ with periodic boundary condition for the spatial domain and zero boundary condition for the velocity domain.

For the Hermite method, we advect not only the distribution function $f$ but also the its gradients of $f$ in $x$ and in $v$ directions.  We have the following equations for $f_x$ and $f_v$,
\begin{equation}
\begin{cases}
\frac{\partial f_x}{\partial t}  +   v  \frac{ \partial f_x }{\partial x} +  \frac{\partial ( E(t,x)f_v ) }{\partial x}= 0,  \\
\frac{ \partial f_v}{\partial t} +  \frac{\partial (vf_x) }{\partial v} +  E(t,x)  \frac{\partial f_v}{\partial v}  =0.
\end{cases}
\label{governing}
\end{equation}

In this section, the SL HWENO scheme solve this system based on the Strang splitting method \cite{cheng1976integration}.
The set of governing equations \eqref{vlasov} and \eqref{governing} of the Strang splitting method is replaced by
\begin{equation}
\label{sl_x}
(SL_x)
\begin{cases}
\frac{\partial f}{\partial t} + v\frac{\partial f}{\partial x} =0    &   (SL_x^0)\\
\frac{\partial f_x}{\partial t} + v\frac{\partial f_x}{\partial x}=0 & (SL_x^1) \\
\frac{\partial f_v}{\partial t} + \frac{ \partial (vf_x)}{\partial v}     =0     & (FD_x)
\end{cases}
\end{equation}
and
\begin{equation}
(SL_v)
\begin{cases}
\frac{\partial f}{\partial t} + E(t,x)\frac{\partial f}{\partial v} =0       & (SL_v^0) \\
\frac{\partial f_v}{\partial t} + E(t,x) \frac{\partial f_v}{\partial v} = 0 &  (SL_v^1)  \\
\frac{\partial f_x}{\partial t} + \frac{\partial( E(t,x)f_v) }{\partial x} =0    &  (FD_v)
\end{cases}.
\end{equation}
 First we solve the system $(SL_x)$ on half a time step, then we solve the system $(SL_v)$ on a time step, and finally solve again the system $(SL_x)$ on half a time step. We will focus our discussion on solving $(SL_x)$ in the following section.

We discretize the computational domain $[a,b]\times[-L,L]$ as
$a=x_{\frac12}<x_{\frac32}<\cdots<x_{N_x+\frac12}=b$,
$-L=v_{\frac12}<v_{\frac32}<\cdots<v_{N_v+\frac12}=L$, with uniformly distributed grid points, i.e.
$x_i=a+(i-\frac12)\Delta x$, $v_j=-L+(j-\frac12)\Delta v$
where grid spacing $\Delta x = x_{i+\frac12}-x_{i-\frac12},\ \Delta v=v_{j+\frac12}-v_{j-\frac12}$.
We let $I_i = [x_{i-\frac12},x_{i+\frac12}] ,\  \forall i=1,\cdots,N_x,\ J_j=[v_{j-\frac12},v_{j+\frac12}] ,\ \forall j=1,\cdots,N_v$ and $T_{ij}= [x_{i-\frac12},x_{i+\frac12}]\times [v_{j-\frac12},v_{j+\frac12}]$.

We let $f_{ij}^n,(f_x)_{ij}^n$ and $(f_v)_{ij}^n$ denote the numerical approximation to the solution $f(x_i,v_j)$, $f_x(x_i,v_j),f_v(x_i,v_j)$ at the time $t^n$ respectively.
Similar to the 1D problem, we introduce $\{\Phi_{i-\frac12,j}^n\}_{ij}$ and $\{\Psi_{i,j-\frac12}^n\}_{ij}$ such that
\begin{equation*}
(f_x)_{ij}^n = \frac{  \Phi_{i+\frac12,j}^n  -\Phi_{i-\frac12,j}^n   }{\Delta x},\ \
(f_v)_{ij}^n = \frac{  \Psi_{i,j+\frac12}^n  -\Psi_{i,j-\frac12}^n   }{\Delta v}
\end{equation*}

%
%
%

In this section, we design a scheme for solving $(SL_x)$ from $t^n$ to $t^{n+1}$. The scheme for $(SL_v)$ would be similar to that for $(SL_x)$.

\noindent
\underline{Initialization:} We use the high order WENO scheme in \cite{jiang1996efficient} to reconstruct $\{\Phi_{i+\frac12,j}^0 \}$ and $\{ \Psi_{i,j+\frac12}^0 \}$ in $x$-direction and $v$-direction respectively.

\noindent
\underline{Update:} We update $\{f_{ij}^{n},\Phi_{i+\frac12,j}^n\}_{ij}$ by the SL HWENO scheme in Section \ref{slhweno1d}.
We update $\Psi_{i,j+\frac12}^n$ by the third equation in \eqref{sl_x} via treating the derivative term as a source term.
In particular,
we apply the following  a central difference scheme coupled with a trapezoid rule for $(FD_x)$,
  \begin{equation}
  \label{fd}
  \begin{split}
  (f_v)_{ij}^{n+1} = (f_v)_{ij}^n - \frac{\Delta t}{2}\left(\frac{ v_{i,j-2}(f_x)_{i,j-2}^n -8 v_{i,j-1}(f_x)_{i,j-1}^n +8v_{i,j+1}(f_x)_{i,j+1}^n -v_{i,j+2}(f_x)_{i,j+2}^n }{24\Delta y}   \right.\\
  \left. \frac{ v_{i,j-2}(f_x)_{i,j-2}^{n+1} -8 v_{i,j-1}(f_x)_{i,j-1}^{n+1} +8v_{i,j+1}(f_x)_{i,j+1}^{n+1} -v_{i,j+2}(f_x)_{i,j+2}^{n+1} }{24\Delta y} \right)
  \end{split}
  \end{equation}
  where
  \begin{equation}
 (f_x)_{ij}^* = \frac{\Phi_{i+\frac12,j}^*-\Phi_{i-\frac12,j}^* }{\Delta x},\ *=n,n+1,
  \end{equation}
and $(f_v)_{ij}^n = \frac{\Psi_{i,j+\frac12}^n - \Psi_{i,j-\frac12}^n }{\Delta y}$.

Equivalently, the scheme \eqref{fd} can be rewritten as  updating $\{\Psi_{i,j-\frac12}^n\}_{ij}$ with
where
    \begin{equation}
    \label{fd1}
  \begin{split}
  \Psi_{i,j-\frac12}^{n+1} =  \Psi_{i,j-\frac12}^n
  -\frac{\Delta t}{2} \left(\frac{-v_{i,j-2}(f_x)_{i,j-2}^n + 7v_{i,j-1}(f_x)_{i,j-1}^n + 7v_{i,j}(f_x)_{i,j}^n -v_{i,j+1}(f_x)_{i,j+1}^n}{24} \right. \\
  \left.  +\frac{-v_{i,j-2}(f_x)_{i,j-2}^{n+1} + 7v_{i,j-1}(f_x)_{i,j-1}^{n+1} + 7v_{ij}(f_x)_{ij}^{n+1} -v_{i,j+1}(f_x)_{i,j+1}^{n+1} }{24}\right).
  \end{split}
  \end{equation}

\begin{remark}
Such source term evolution by the central difference and a trapezoid rule for time integration  has time step restriction and may cause instability if the time step size is too large.
\end{remark}
 \begin{remark}
There is another form of governing equations for $f$, $f_x$ and $f_v$ in \cite{besse2003semi,besse2008convergence},
\begin{equation}
\begin{cases}
\frac{\partial f}{\partial t} + v \frac{\partial f}{\partial x} + E(t,x)\frac{\partial f}{\partial v}=0,  \\
\frac{\partial f_x}{\partial t}  +   v  \frac{ \partial f_x }{\partial x} + \frac{\partial E(t,x)}{\partial x}f_v + E(t,x)\frac{\partial f_x}{\partial v}= 0,  \\
\frac{ \partial f_v}{\partial t}+ f_x + v \frac{\partial f_v}{\partial x}  +  E(t,x)  \frac{\partial f_v}{\partial v}  =0.
\end{cases}
\label{governing1}
\end{equation}
Then in the context of operator splitting,  the third equation in \eqref{sl_x} will be
\begin{equation}
\frac{ \partial f_v}{\partial t}+ f_x + v \frac{\partial f_v}{\partial x}  = 0.
\end{equation}
We observe that we can only design the SL scheme on $\frac{ \partial f_v}{\partial t}+ v \frac{\partial f_v}{\partial x}=0$.
It is not mass conservative for $f_v$ with the source term $f_x$. It lead to  the difficulty for writing $(f_v)_{ij}^n$ in a flux difference form as
$(f_v)_{ij}^n = \frac{  \Psi_{i,j+\frac12}^n  -\Psi_{i,j-\frac12}^n   }{\Delta v}$

\end{remark}

This SL HWENO scheme may lead to oscillations in large time stepping size when it simulates the VP system mainly from the source term evolution \eqref{fd1}.
We propose to apply the WENO limiter \cite{qiu2004hermite,qiu2005runge}
before HWENO evolution as a pre-processing procedure which is similar to the procedure in \cite{guo2013hybrid}. We use the TVB limiter \cite{cockburn2001runge,qiu2004hermite,qiu2005runge}
with problem dependent TVB constants to identify troubled cells.
For details of the procedure of the limiter, we refer to \cite{qiu2005runge}. Below we provide the flow chart of the conservative SL HWENO with WENO limiters for the VP simulations.

\begin{algorithm}
Conservative SL HWENO scheme for the VP system.

\begin{description}
  \item [Step 1.] Apply WENO limiter as a pre-processing procedure to reconstruct $\Phi_{i+\frac12,j}^n$.

  First, we use TVB limiter to identify the ¡°troubled cells,¡± namely, those cells which might
need the limiting procedure. Let:
\begin{equation}
\tilde{f}_{ij} = \Phi_{i+\frac12,j}^n - f_{ij}^n,\  \ \tilde{\tilde{f}}_{ij}=-\Phi_{i-\frac12,j}^n + f_{ij}^n
\end{equation}
These are modified by the modified minmod function;
\begin{equation}
\tilde{f}_{ij}^{(mod)} = \tilde{m}\left( \tilde{f}_{ij} , f_{i+1,j}^n - f_{ij}^n , f_{ij}^n -f_{i-1,j}^n \right),
\end{equation}
\begin{equation}
\tilde{\tilde{f}}_{ij}^{(mod)} = \tilde{m}\left( \tilde{\tilde{f}}_{ij} , f_{i+1,j}^n - f_{ij}^n , f_{ij}^n -f_{i-1,j}^n \right),
\end{equation}
where $\tilde{m}$ is given by
\begin{equation}
\tilde{m}(a_1,a_2,\cdots, a_n)=
\begin{cases}
a_1 &\text{if} \ |a_1|\leq M_x(\Delta x)^2, \\
m(a_1,a_2,\cdots,a_n) &  otherwise,
\end{cases}
\end{equation}
and the minmod function $m$ is given by
\begin{equation}
m(a_1,a_2,\cdots,a_n)=
\begin{cases}
s\cdot \min_{1\leq j \leq n} |a_j| & \text{if}\ sign(a_1)=sign(a_2)=\cdots=sign(a_n)=s,\\
0                                  &  otherwise.
\end{cases}
\end{equation}
The TVB limiter parameter $M_x>0$ is a constant.
If $\tilde{f}_{ij}^{(mod)} \neq \tilde{f}_{ij}$ or $ \tilde{\tilde{f}}_{ij}^{(mod)} \neq \tilde{\tilde{f}}_{ij} $, we declare the $T_{ij}$ as a troubled cell.

Then we replace $\Phi_{i+\frac12,j}^n$ and $\Phi_{i-\frac12,j}^n$ in those troubled cells by
WENO reconstruction.

  \item [Step 2.] Perform a half time step advection in physical space.
  \item [Step 3.] Compute the electric field at the half step by substituting $f^*$ into equation \eqref{eletric} and solve for $E^*(x)$.
  \item [Step 4.] Similar to Step 1 , we apply WENO limiter as a pre-processing procedure to reconstruct $\Psi_{i,j+\frac12}^*$.

  \item [Step 5.] Perform a full time step advection in velocity space.
  \item [Step 6.] We use the same pre-processing procedure like Step 1 to reconstruct $\Phi_{i+\frac12,j}^{**}$.
  \item [Step 7.] Perform a half time step advection in physical space.
\end{description}

\end{algorithm}

\section{Numerical tests}
\label{tests}
In Section \ref{test1}, we first test the performances of the proposed SL HWENO schemes for the 1D transport problem, and then the rigid body rotation was tested by the Strang splitting conservative SL HWENO scheme. In Section \ref{testvp}, we demonstrate
the utility of the SL HWENO scheme by applying it to classical problems from plasma physics, such
as Landau damping and two-steam instability.

\subsection{Test problems}
\label{test1}

\begin{example}
1D transport problem.

Consider the linear advection equation:
\begin{equation}
\label{linear}
f_t+f_x = 0,\   x\in[0,2\pi].
\end{equation}
The conservative SL methods with fifth order HWENO reconstruction is used to solve equation \eqref{linear}.
Table \ref{table:22} gives the $L_1$ error, and the corresponding order of convergence, of both the SL HWENO scheme and the SL WENO scheme when applied to equation \eqref{linear} with smooth initial data
\begin{equation}
f(x,0)=\sin(x) \  \ (0\leq x\leq 2\pi).
\end{equation}
We observe that these two methods have fifth order accuracy, and the SL HWENO scheme is more accurate than the SL WENO scheme.
\begin{table}[!ht]\small
\caption{\small Order of accuracy for \eqref{linear} with $f(x,t=0)=\sin(x)$ at $T=20. \ CFL =1.2.$ }
\vspace{0.1in}
\centering
\begin{tabular}{|l|ll|ll|}
\hline
{\multirow{2}{*}{N}} & \multicolumn{2}{c|}{HWENO}  &\multicolumn{2}{c|}{WENO}   \\         
\multicolumn{1}{|l|}{}&
 \multicolumn{1}{|l}{$L_1$ error} & Order  & {$L_1$ error} & Order \\
\hline
  32  &     4.03E-05 & -           &       7.31E-05  &  -     \\ \hline
  64  &     1.17E-06 &     5.10    &       2.23E-06  &  5.03  \\ \hline
  96  &     1.52E-07 &     5.05    &       2.93E-07  &  5.00  \\ \hline
 128  &     3.56E-08 &     5.04    &       6.97E-08  &  5.00  \\ \hline
 160  &     1.16E-08 &     5.03    &       2.28E-08  &  5.00  \\ \hline
 192  &     4.62E-09 &     5.03    &       9.16E-09  &  5.00  \\ \hline
\end{tabular}
\label{table:22}
\end{table}

Next, to evaluate the capability of the schemes in capturing both discontinuity and smooth solution, we test \eqref{linear} with the initial distribution including four types of profiles
\begin{equation}
f(x,0)=
\begin{cases}
\frac16 (G(x,\beta,z-\delta) + G(x,\beta,z+\delta) + 4G(x,\beta,z)   ) & \text{for} \  -0.8\leq x\leq-0.6, \\
1                               &  \text{for}\  -0.4\leq x \leq -0.2,\\
1-| 10(x-0.1) |                 &  \text{for}\  0.0\leq x \leq 0.2,\\
\frac16( F(x,\alpha,a-\delta) + F(x,\alpha,a+\delta) + 4F(x,\alpha,a) ) & \text{for} \ 0.4\leq x\leq 0.6, \\
0                               &     \text{otherwise},
\end{cases}
\end{equation}
where $G(x,\beta,z)=e^{-\beta(x-z)^2}$ and $F(x,\alpha,a)=\sqrt{\max(1-\alpha^2(x-a)^2,0) }$.
The constants are specified as $a=0.5, z= -0.7,\delta=0.005,\alpha=10$ and $\beta=\frac{\log2}{36\delta^2}$. The boundary condition is periodic.
We compute the solution up to $t=8$ with 200 points. The results are shown in Figure \ref{Fig:eg1}. Non-oscillatory numerical capture of discontinuities is observed.
\begin{figure}[h!]
\centering                              
\includegraphics[height=80mm]{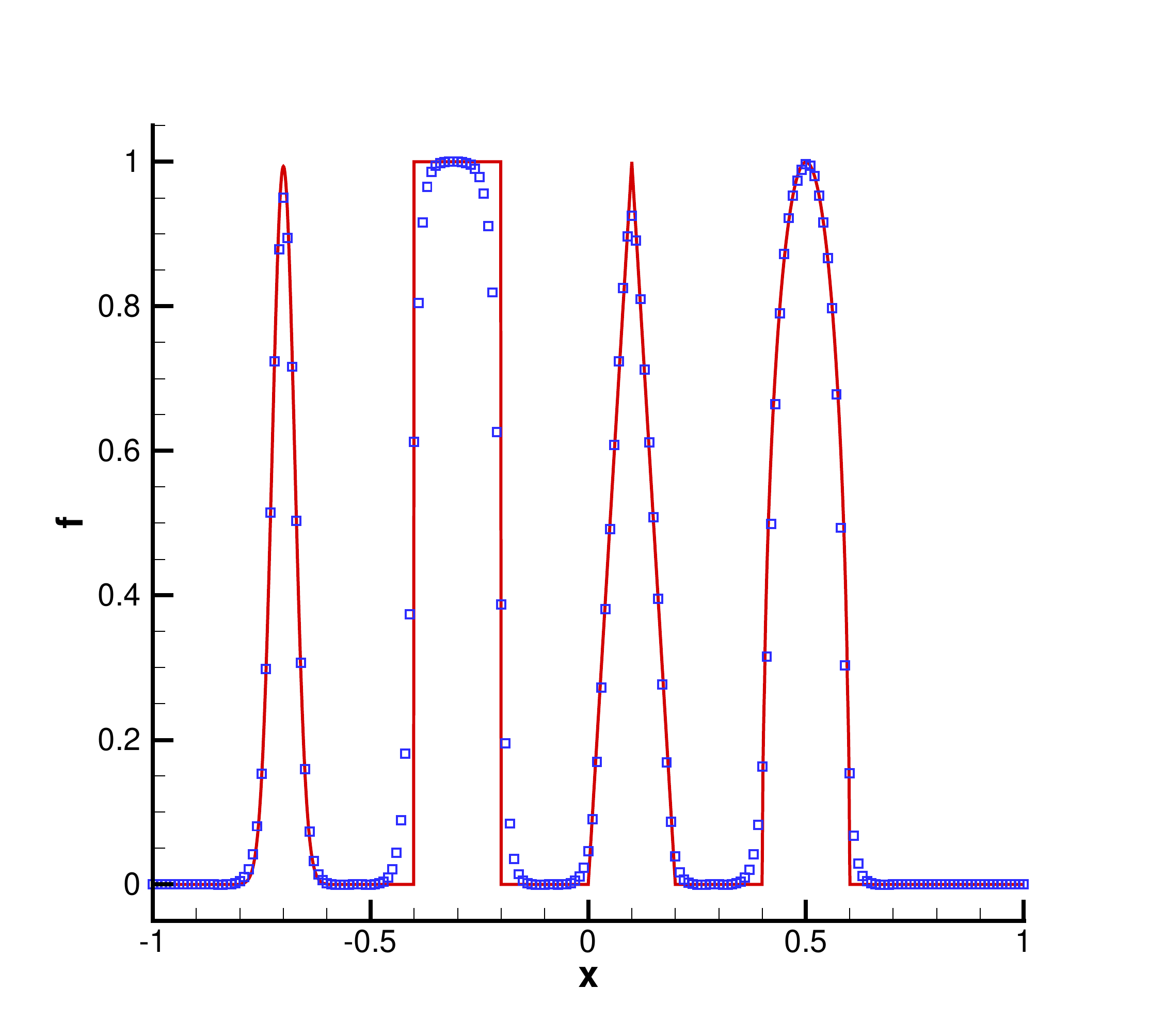}
\caption{\small The SL HWENO scheme; 200 points, CFL=1.2, T=8. }
\label{Fig:eg1}
\end{figure}
\end{example}

\begin{example}
Rigid body rotation.

Consider the rigid body rotation
\begin{equation}
f_t + v_x(x,y) f_x + v_y(x,y) f_y = 0.
\label{rotation}
\end{equation}

First, we consider a smooth case for accuracy test.  The domain is $[-2\pi,2\pi]\times[-2\pi,2\pi]$. The velocity field is given by $v_x(x,y)=-y,\ v_y(x,y)=x$.
Table \ref{table:11} gives the error and convergence rates of the scheme for the time step $\Delta t =CFL/(\frac{2\pi}{\Delta x}+\frac{2\pi}{\Delta y} )$ with $CFL=1.2$ for smooth initial data $f(x,y,0)=\exp(-x^2-y^2)$. The high order convergence of the scheme  is observed.

\begin{table}[!ht]\small
\caption{\small Order of accuracy for \eqref{rotation} with $f(x,y,t=0)=\exp(-x^2-y^2)$ at $T=2\pi. \ CFL=1.2$. }
\vspace{0.1in}
\centering
\begin{tabular}{|l|ll|ll|ll|}
\hline
\multicolumn{1}{|l|}{ N }&
 \multicolumn{1}{|l}{$L_1$ error} & Order  & {$L_2$ error} & Order & {$L_\infty$ error} & Order\\
\hline
      20 &     1.31E-02 & &     2.49E-02 & &     1.19E-01 & \\ \hline
    40 &     1.05E-03 &     3.65 &     1.90E-03 &     3.71 &     1.50E-02 &     2.99 \\ \hline
       80 &     4.34E-05 &     4.59 &     8.88E-05 &     4.42 &     4.99E-04 &     4.91 \\ \hline
     160 &     2.03E-06 &     4.42 &     3.97E-06 &     4.48 &     2.11E-05 &     4.56 \\ \hline
     320 &     6.50E-08 &     4.96 &     1.36E-07 &     4.87 &     7.33E-07 &     4.85 \\ \hline
\end{tabular}
\label{table:11}
\end{table}

Secondly, we consider a test case introduced in \cite{leveque1996high}. The domain is $[-0.5,0.5]^2$. The velocity field is given by
\begin{equation}
v_x(x,y) = -2\pi y,\   \ v_y(x,y)=2\pi x.
\label{velocity}
\end{equation}
The initial condition we used is plotted in Figure \ref{Fig:11}. It includes a slotted disk, a cone and a smooth hump. The numerical solutions for the time step $\Delta t =CFL/(\frac{\pi}{\Delta x}+\frac{\pi}{\Delta y} )$ with $CFL=1.2$ after one full revolutions by the conservative SL HWENO scheme (denoted as CSLHWENO-WO) are plotted in Figure \ref{Fig:11}.
 However, the solution of the scheme without WENO limiter with $CFL=2.2$ will be oscillatory. We apply the WENO limiter to the scheme (denoted as CSLHWENO-WL). The numerical solution and trouble cells at the last time step are presented in Fig. \ref{Fig:rotaion1d4}. In Fig. \ref{Fig:13}, we plot the 1D cut of the solution compared with the exact solution.
Non-oscillatory capturing of discontinuities is observed.

\begin{figure}[h!]
\centering
\subfigure[The initial profile]{
\includegraphics[height=70mm]{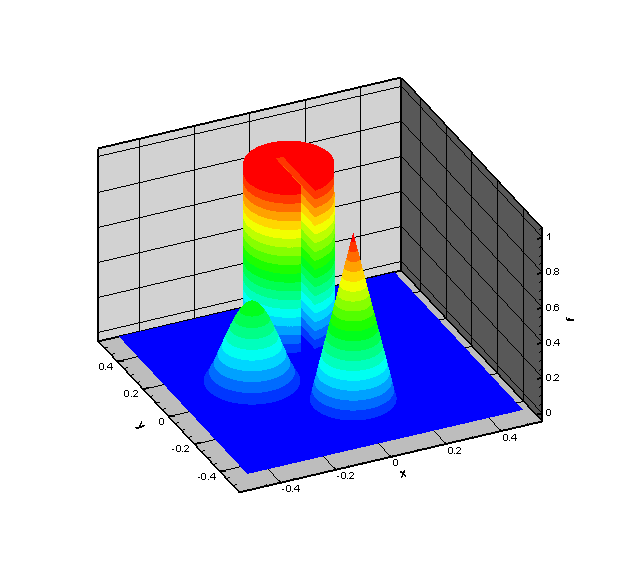}  }
\subfigure[SLHWENO with $CFL=1.2$]{
\includegraphics[height=70mm]{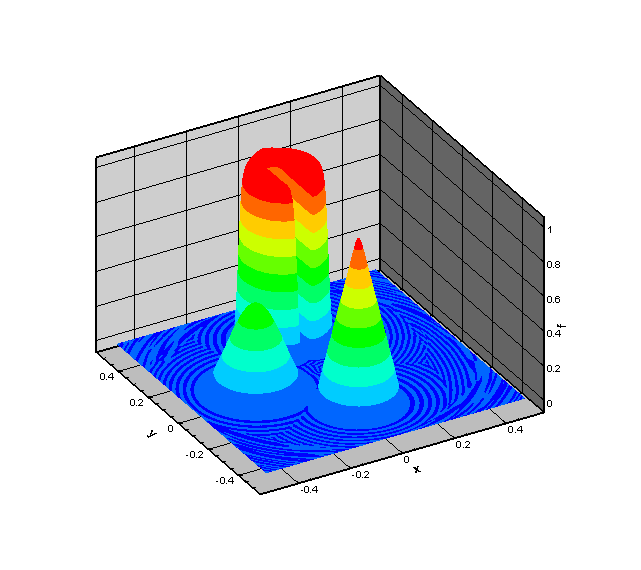} }
\caption{Left: Plots of the initial profile. Right: Plots of the numerical solution for equation \eqref{rotation} with the velocity field \eqref{velocity}; $CFL=1.2$;  $T=1$;  The numerical mesh has a resolution of $200\times200$; Conservative SL HWENO scheme without WENO limiter.}
\label{Fig:11}
\end{figure}

\begin{figure}[h!]
\centering
\includegraphics[height=70mm]{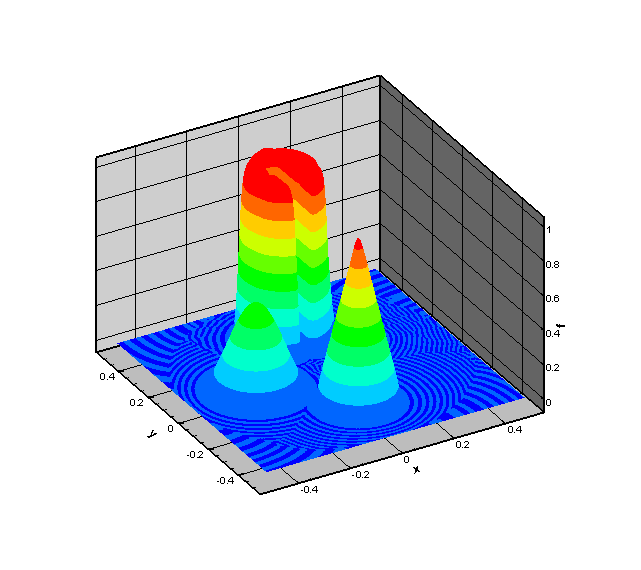}
\includegraphics[height=70mm]{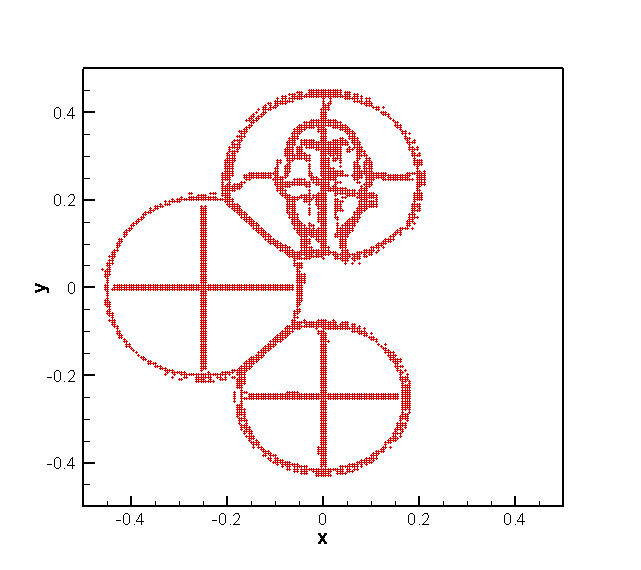}
\caption{Left: The numerical solution for equation \eqref{rotation} with the velocity field \eqref{velocity}.
Right: Trouble cells. $CFL=2.2$. TVB constant $M =1.0$. $T=1.$ The numerical mesh has a resolution of $200\times200$. Conservative SL HWENO scheme with WENO limiter.}
\label{Fig:rotaion1d4}
\end{figure}

\begin{figure}[h!]
\centering
\subfigure[y=-0.25]{
\includegraphics[height=65mm]{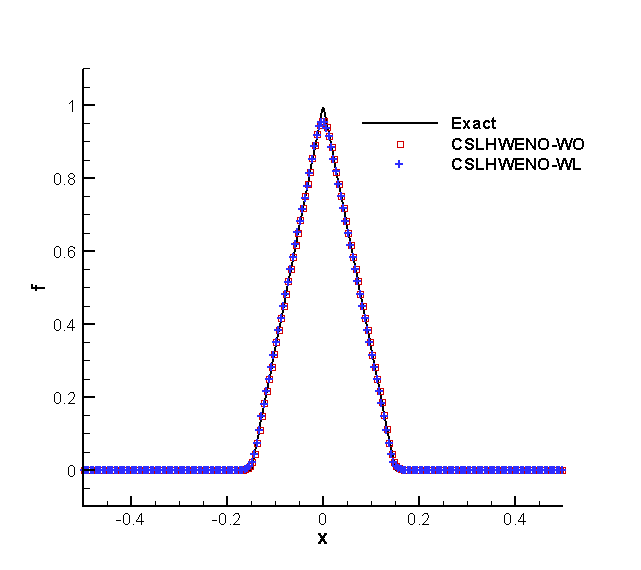} }
\subfigure[y=0.25]{
\includegraphics[height=65mm]{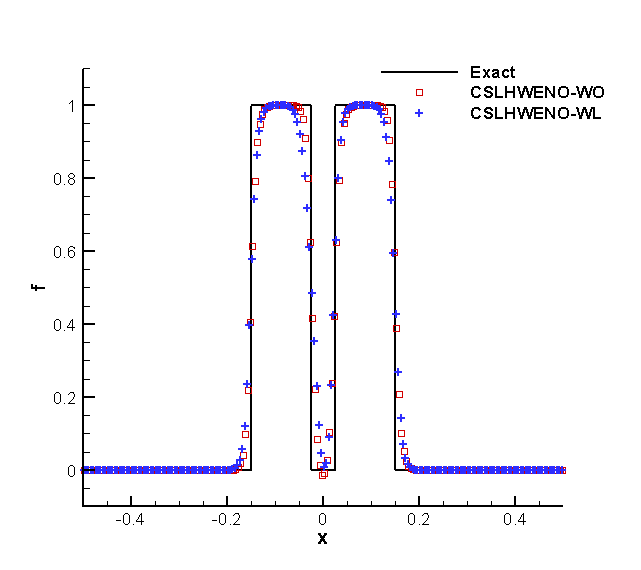} }
\subfigure[x=-0.25]{
\includegraphics[height=65mm]{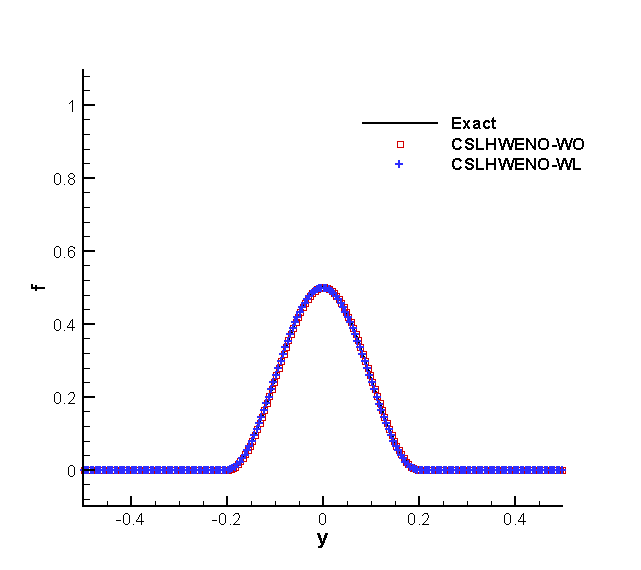} }
\subfigure[x=0]{
\includegraphics[height=65mm]{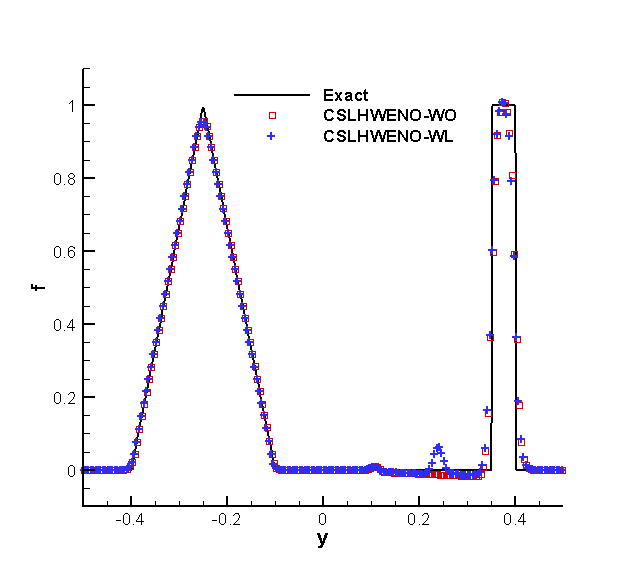} }
\caption{Plots of the 1D cuts of the numerical solution for equation \eqref{rotation} at $y=-0.25,y=0.25,x=-0.25,x=0$ (from top left to bottom right). The solid line depicts the exact solution. The numerical mesh has a resolution of $200\times200$.}
\label{Fig:13}
\end{figure}

\end{example}

\subsection{The VP system}
\label{testvp}
In this subsection, we apply the conservative SL HWENO scheme  to the VP system.
Periodic boundary conditions are imposed in the $x$-direction and zero boundary conditions are imposed in the $v$-direction for all of our test problems. Because of the periodicity in space, a fast Fourier transform (FFT) is used to solve the 1-D Poisson equation. $\rho(x,t)$ is computed by the rectangular rule, $\rho(x,t)=\int f(x,v,t)dv = \mathop{\sum_j}f(x,v_j,t)\Delta v$, which is spectrally accurate \cite{boyd2001chebyshev}, when the underlying function is smooth enough.
We recall several norms in the VP system below, which remain constant in time.
\begin{enumerate}
  \item  $L^p$ norm $1\leq p \leq \infty$:
           \begin{equation}
           \| f \|_p = \left( \int_v \int_x |f(x,v,t)|^p dxdv  \right)^{\frac{1}{p} }
           \end{equation}

  \item Energy:
                \begin{equation}
                Energy = \int_v\int_x f(x,v,t)v^2 dxdv + \int_x E^2(x,t) dx,
                \end{equation}
                where $E(x,t)$ is the electric field.
  \item Entropy:
               \begin{equation}
               Entropy = \int_v\int_x f(x,v,t) \log ( f(x,v,t) ) dxdv.
               \end{equation}
\end{enumerate}
Tracking relative deviations of these quantities numerically will be a good measurement of the quality of numerical schemes.
The relative deviation is defined to be the deviation away from the corresponding initial value divided by the magnitude of the
initial value. It is expected that our scheme will conserve mass.
However, the positivity of $f$ will not be preserved. Thus, when numerically computing the entropy, we compute
$\int_v\int_x f(x,v,t) \log | f(x,v,t) | dxdv$. We set the time step by $\Delta t =CFL /( v_{max}/\Delta x +\max(E(x))/\Delta y )$, where $v_{max}$ is the maximum velocity on the phase space mesh.

In the following, we test the conservative SL HWENO scheme with $CFL=1.2$, denoted as "CSLHWENO-WO", to solve the VP system. This schemes will be compared with the fifth order conservative SL WENO scheme proposed in \cite{qiu2010conservative} with the same $CFL=1.2$, denoted as "CSLWENO".
Moreover, we will study the conservative SL HWENO scheme with the large $CFL=2.2$. In this case, the WENO limiter in certain TVB constants is needed to enforce the stability of this scheme, and denote the scheme as "CSLHWENO-WL".

\begin{example}
Weak Landau damping.

Consider the weak Landau damping for the VP system. The initial condition used here is,
\begin{equation}
f(x,v,t=0) = \frac{1}{\sqrt{2\pi} } (1+\alpha \cos(kx) )\exp \left(-\frac{v^2}{2} \right),
\label{Landau}
\end{equation}
with $\alpha=0.01$ and $k=0.5$. Our simulation parameters are $v_{max}=5,\ N_x = 64,\ N_v =128$.  The time evolution of the $L^2$ and $L^\infty$ norms of the electric field (in semi-log scale) are plotted in the upper plots of Figure \ref{weakLandau}. The correct damping of the electric field of CSLWENO and CSLHWENO-WO is observed in the plots, benchmarked with the  theoretical value $\gamma=0.1533$ \cite{filbet2001conservative} (the solid line in the same plots). We observe that the conservative SL scheme generates very consistent results, performing very well in recovering the damping rate. Time evolution of the relative deviations of  the $L^1$, $L^2$ solution norms, energy, entropy in the discrete sense are demonstrated in the middle and bottom plots in Figure \ref{weakLandau}. The advantage of using conservative schemes in preserving the relevant physical norms is observed.
CSLHWENO-WO is observed to perform slightly better than the CSLWENO in preserving norms.

\begin{figure}[h!]
\centering
\includegraphics[height=65mm]{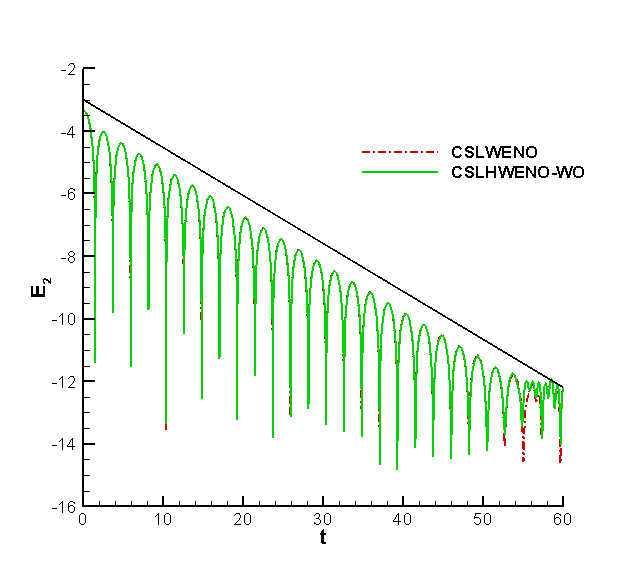}
\includegraphics[height=65mm]{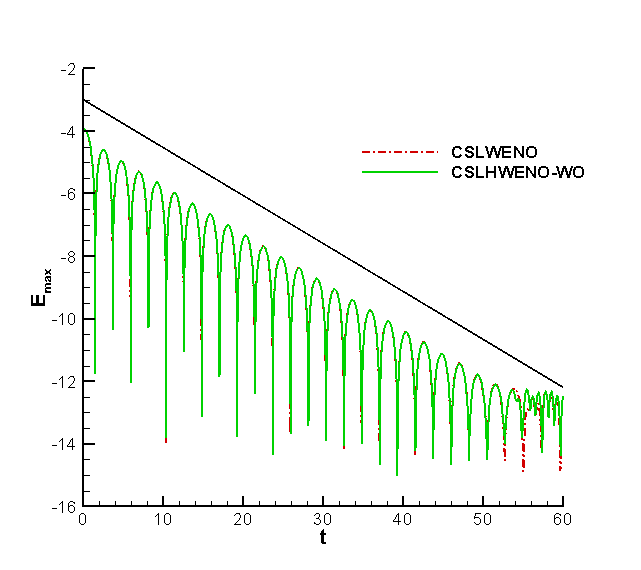}
\includegraphics[height=65mm]{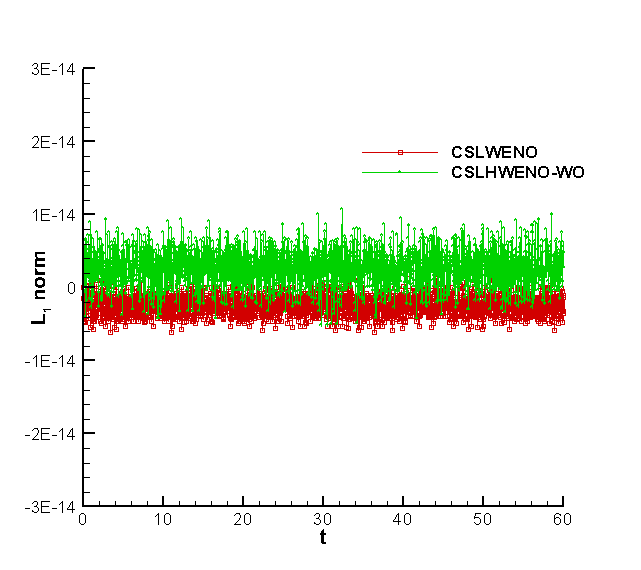}
\includegraphics[height=65mm]{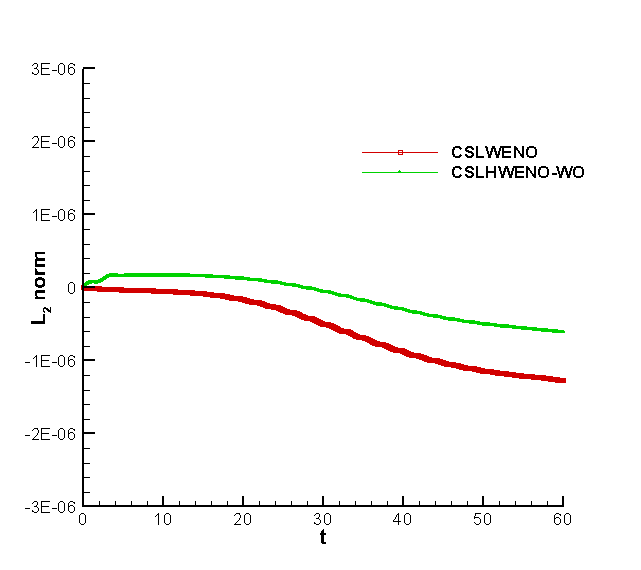}
\includegraphics[height=65mm]{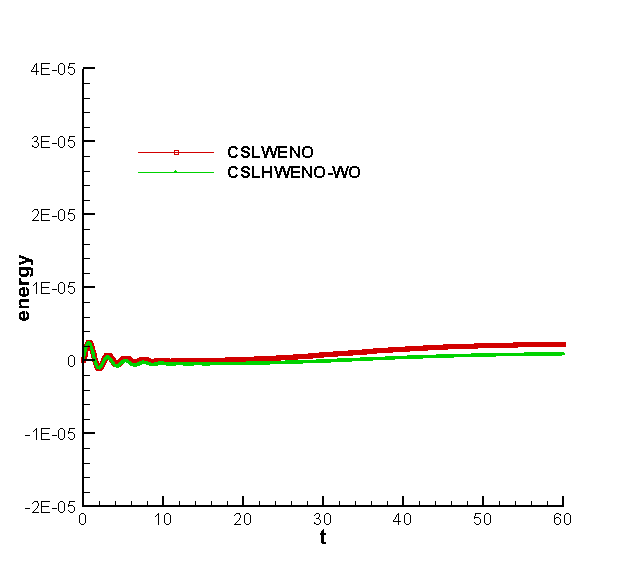}
\includegraphics[height=65mm]{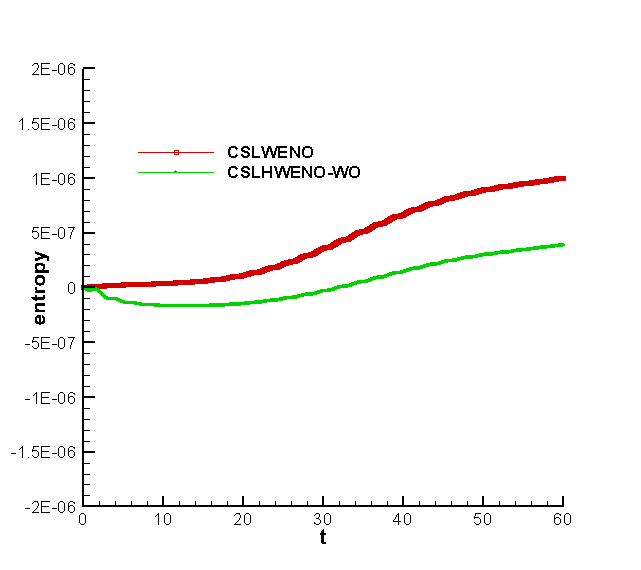}
\caption{Weak Landau damping: time evolution of the electric field in $L^2$ (upper left) and $L^\infty$ (upper right) norms, Time evolution of the relative deviations of $L^1$ (middle left) and $L^2$ (middle right) norms of the solution as well as the discrete kinetic energy (lower left) and entropy (lower right).}
\label{weakLandau}
\end{figure}

\end{example}

\begin{example}
Strong Landau damping.

Consider the strong Landau damping for the VP system. We simulate the VP system with the initial condition in equation \eqref{Landau} with $\alpha=0.5$ and $k=0.5$. Our simulation parameters are $v_{max}=5,\ N_x = 128,\ N_v =256$.
In the first row of Figure \ref{strongLandau}, the time evolution of the $L^2$ and $L^\infty$ norms of the electric field with the linear decay rate $\gamma_1=-0.2812$ and $\gamma_2=0.0770$ \cite{cheng1976integration,guo2013hybrid}  (in semi-log scale) are plotted.
 The Time evolution of the relative deviations of discrete $L^1$ norm, $L^2$ norm, kinetic energy and entropy for CSLWENO and CSLHWENO-WO are plotted in the second and third rows of Figure \ref{strongLandau}.
 CSLHWENO-WO scheme is observed to perform slightly better in preserving these norms than the CSLWENO scheme.
 In Figure \ref{strongLandau1}, numerical solutions of CSLWENO  and CSLHWENO-WO  at $t=30$ are plotted. Compared with CSLWENO, slightly  better resolution is observed for CSLHWENO-WO.
 The results of conservative SL scheme with $CFL=2.2$ by WENO limiter as well as the trouble cells in the last step of the CSLHWENO-WL evolution are presented in the bottom plots of Figure \ref{strongLandau1}.
 The results of CSLHWENO-WL ($CFL=2.2$) are observed to be  comparable to those of CSLHWENO-WO ($CFL=1.2$).
It is also observed that when the filamentation structures become under-resolved by the numerical mesh, trouble cells are identified and WENO limiters are applied.

\begin{figure}[h!]
\centering
\includegraphics[height=65mm]{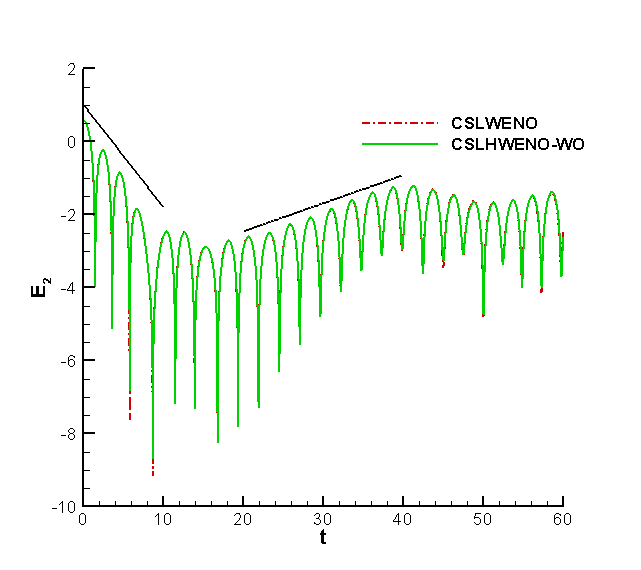}
\includegraphics[height=65mm]{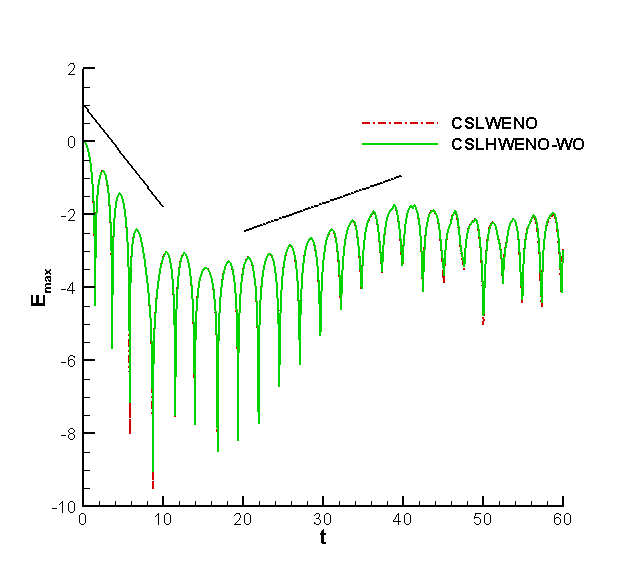}
\includegraphics[height=65mm]{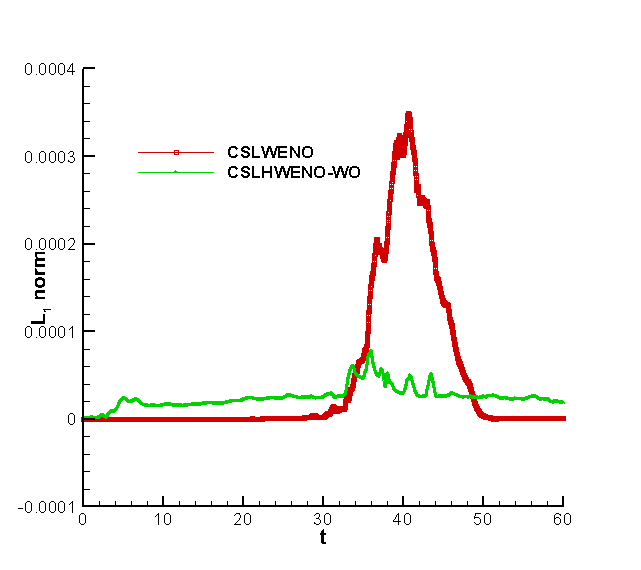}
\includegraphics[height=65mm]{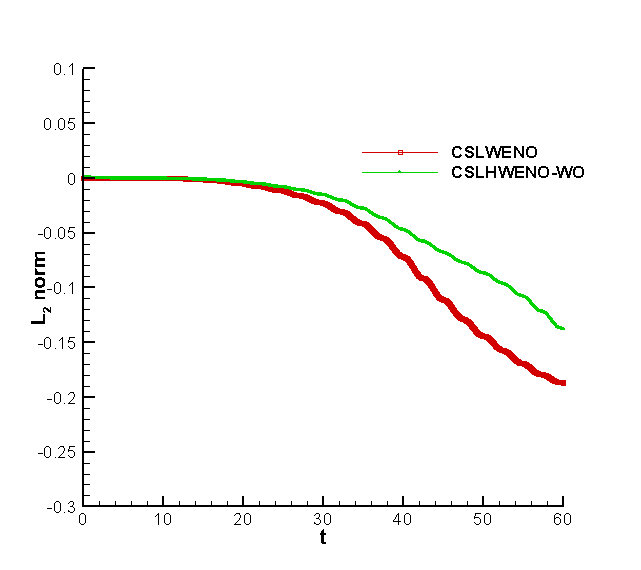}
\includegraphics[height=65mm]{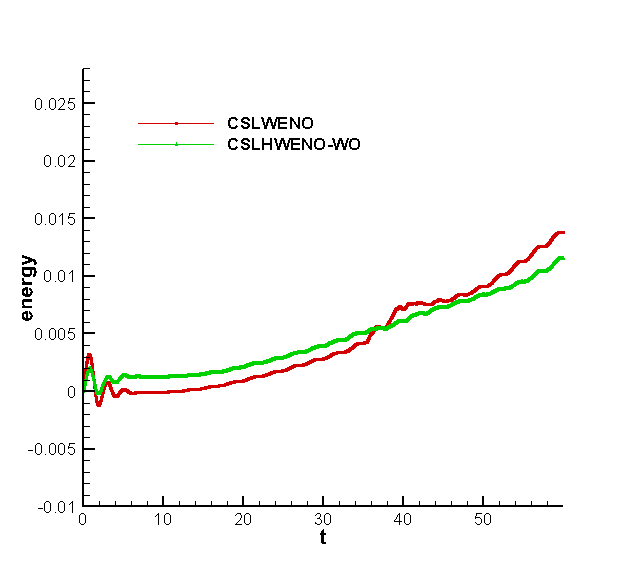}
\includegraphics[height=65mm]{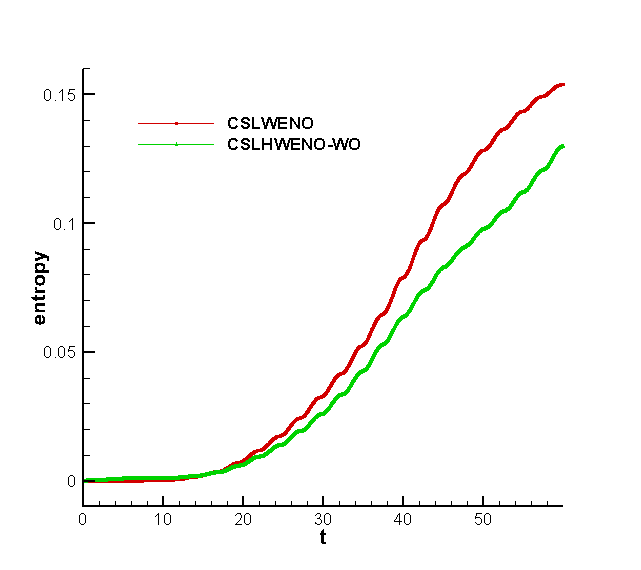}
\caption{Strong Landau damping: time evolution of the electric field in $L^2$ (upper left) and $L^\infty$ (upper right) norms, time evolution of the relative deviations of $L^1$ (middle left) and $L^2$ (middle right) norms of the solution as well as the discrete kinetic energy (lower left) and entropy (lower right).}
\label{strongLandau}
\end{figure}
\begin{figure}[h!]
\centering
\subfigure[CSLWENO]{
\includegraphics[height=70mm]{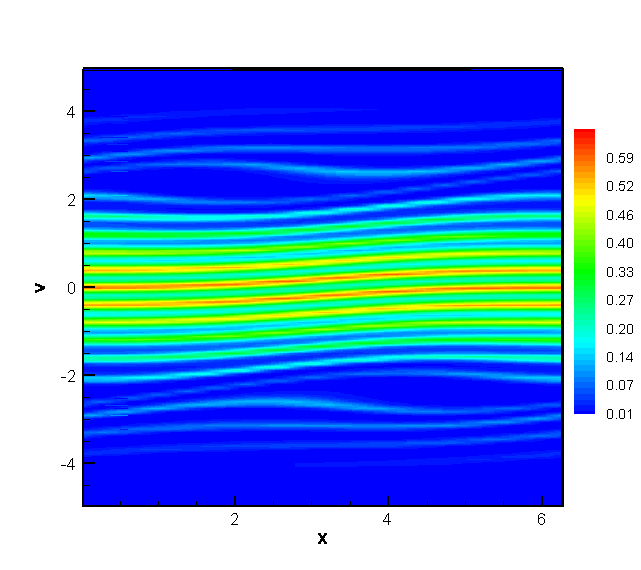}  }
\subfigure[CSLHWENO-WO]{
\includegraphics[height=70mm]{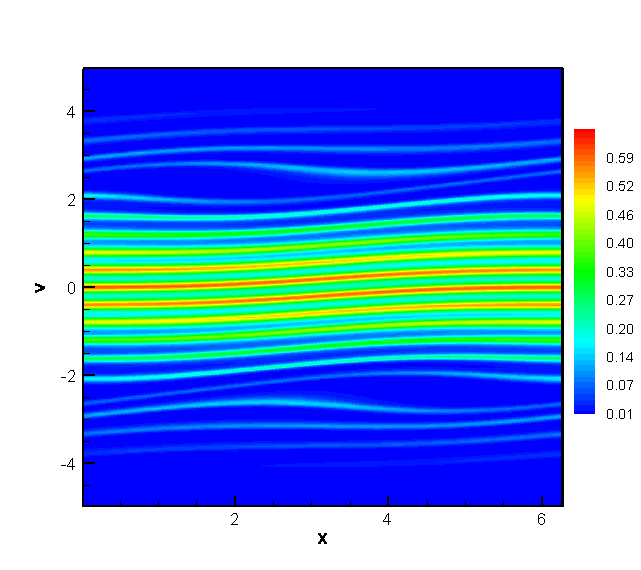} }
\subfigure[CSLHWENO-WL]{
\includegraphics[height=70mm]{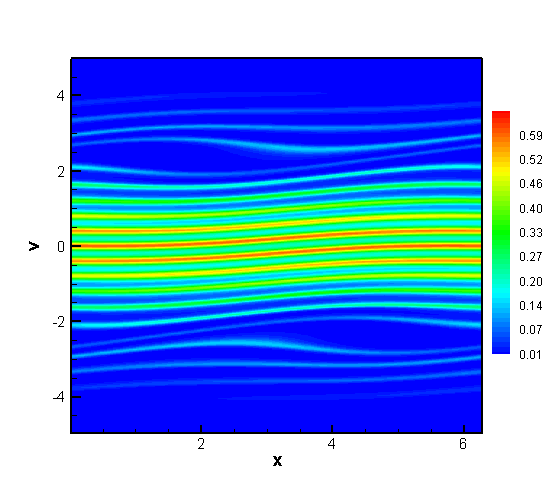} }
\subfigure[trouble cells]{
\includegraphics[height=70mm]{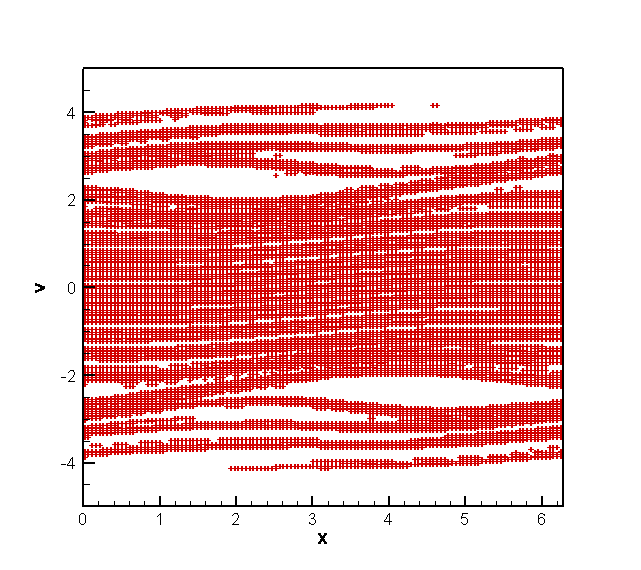} }
\caption{Strong Landau damping. $T=30$. $N_x \times N_v = 128\times256$. Top left: CSLWENO. Top right: CSLHWENO-WO ($CFL=1.2$). Bottom left: CSLHWENO-WL ($CFL=2.2$); the TVB constants $M_x=M_y=1$. Bottom right: trouble cells of CSLHWENO-WL at the last time step.}
\label{strongLandau1}
\end{figure}

\end{example}

\begin{example}
Two stream instability \cite{filbet2003comparison}.

Consider two stream instability, with an unstable initial distribution function,
\begin{equation*}
f(x,v,t=0) = \frac{2}{7\sqrt{2\pi} } (1+5v^2) ( 1 + \alpha \left(  ( \cos(2kx)+\cos(3kx) )/1.2 + \cos(kx) \right) )\exp\left(-\frac{v^2}{2} \right)
\end{equation*}
with $\alpha =0.01,\ k=0.5$. The length of the domain in the $x$ direction is $L=\frac{2\pi}{k}$ and the background ion distribution function is fixed, uniform and chosen so that the total net charge density for the system is zero. Our numerical simulation parameters are $v_{max}=5,\ N_x =64,\ N_v=128$.
 In the first row of Figure \ref{two1}, the time evolution of the $L^2$ and $L^\infty$ norms of the electric field (in semi-log scale) for CSLWENO and CSLHWENO-WO are presented. Comparable results of the discrete $L^1$ norm, $L^2$ norm, kinetic energy and entropy for CSLWENO and CSLHWENO-WO are observed  in the second and third rows of Figure \ref{two1}.
In Figure \ref{two}, we show  numerical solutions at $T=53$ of CSLWENO, CSLHWENO-WO
  and CSLHWENO-WL.
  Comparable solutions are observed for the CSLWENO ($CFL=1.2$), CSLHWENO-WO ($CFL=1.2$) and CSLHWENO-WL ($CFL=2.2$). Note that the solution of CSLHWENO-WO at $CFL=2.2$ (not presented here) will be oscillatory without the limiter.
\begin{figure}[h!]
\centering
\includegraphics[height=65mm]{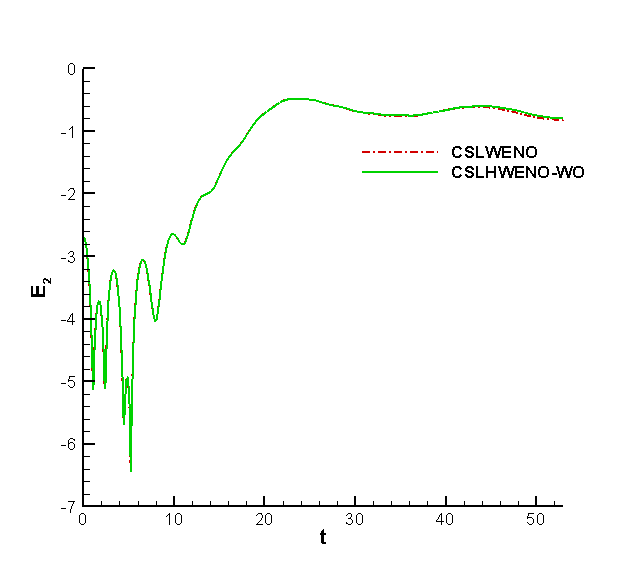}
\includegraphics[height=65mm]{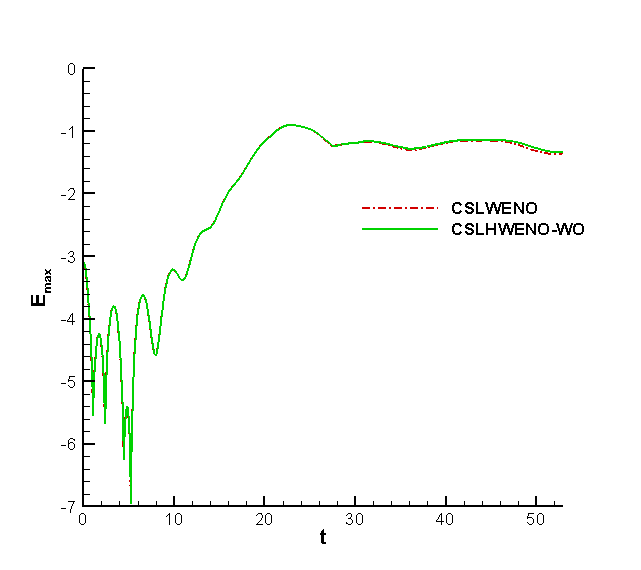}
\includegraphics[height=65mm]{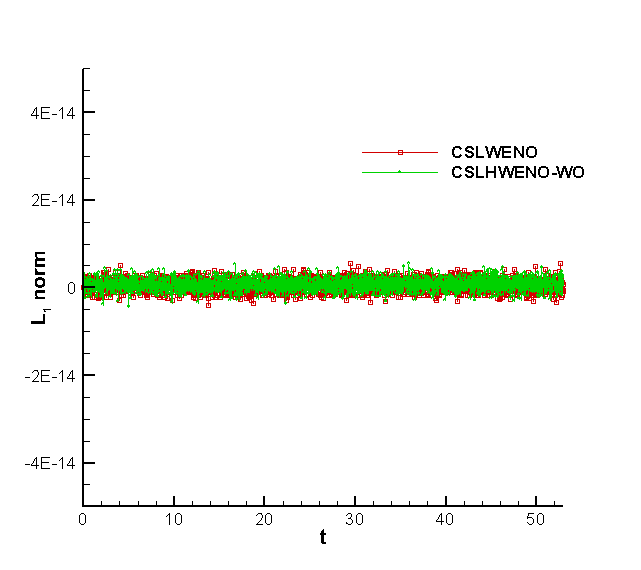}
\includegraphics[height=65mm]{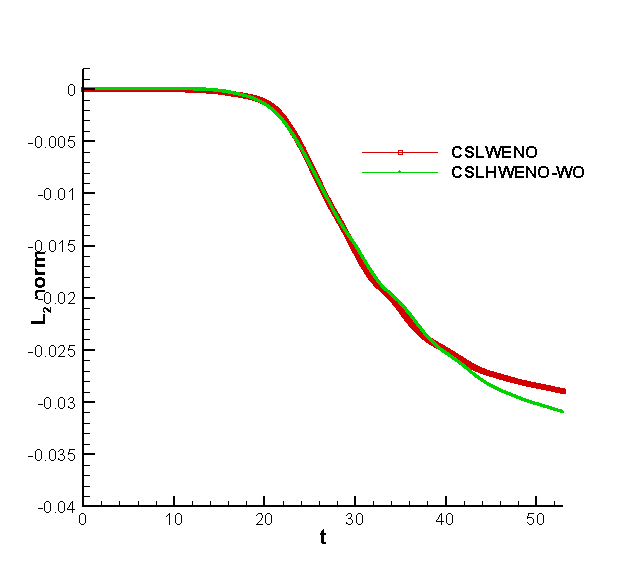}
\includegraphics[height=65mm]{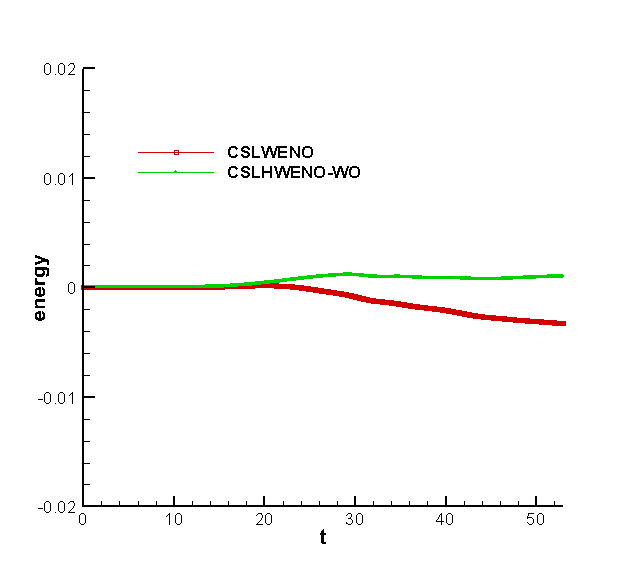}
\includegraphics[height=65mm]{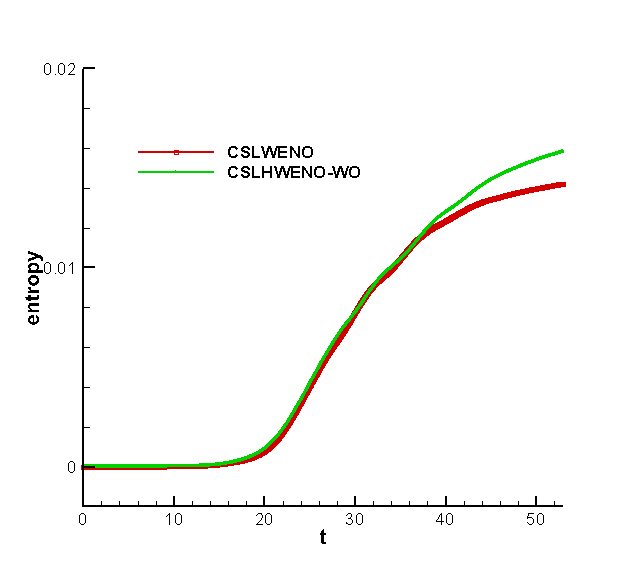}
\caption{Two-stream instability: time evolution of the electric field in $L^2$ (upper left) and $L^\infty$ (upper right) norms, time evolution of the relative deviations of $L^1$ (middle left) and $L^2$ (middle right) norms of the solution as well as the discrete kinetic energy (lower left) and entropy (lower right) for CSLWENO and CSLHWENO-WO.}
\label{two1}
\end{figure}

\begin{figure}[h!]
\centering
\subfigure[CSLWENO]{
\includegraphics[height=70mm]{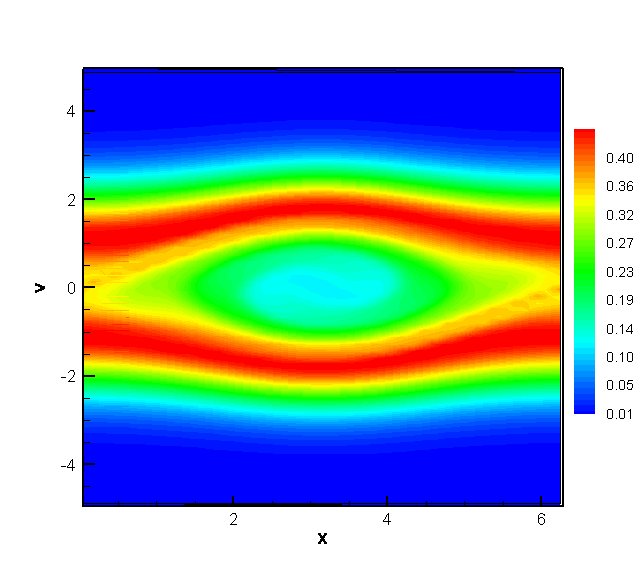} }
\subfigure[CSLHWENO-WO]{
\includegraphics[height=70mm]{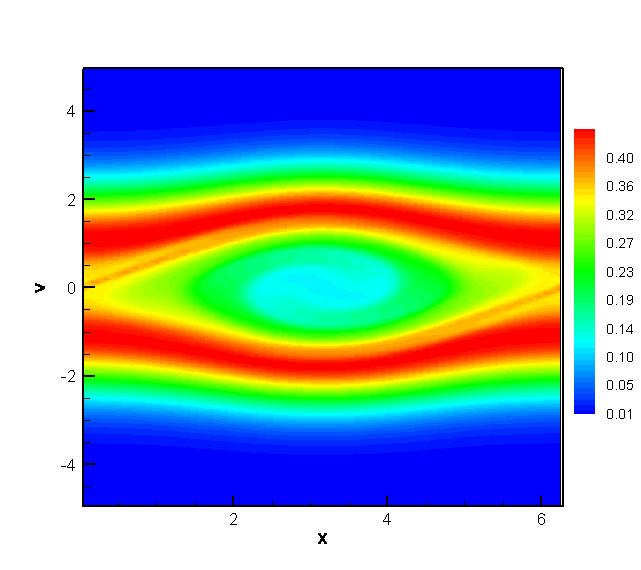} }
\subfigure[CSLHWENO-WL]{
\includegraphics[height=70mm]{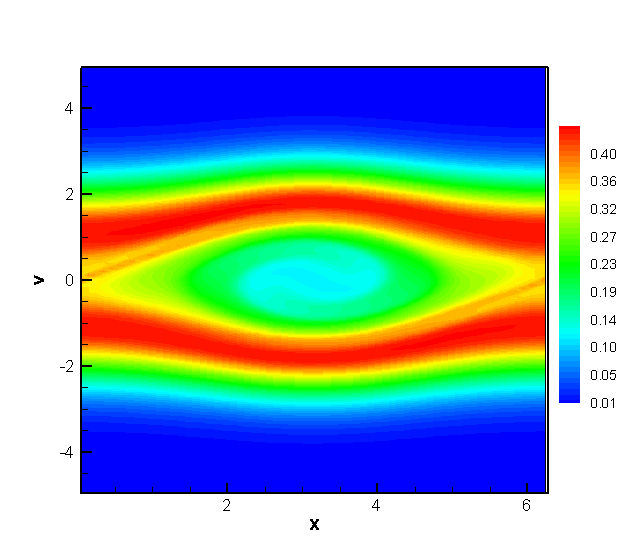} }
\subfigure[trouble cells]{
\includegraphics[height=70mm]{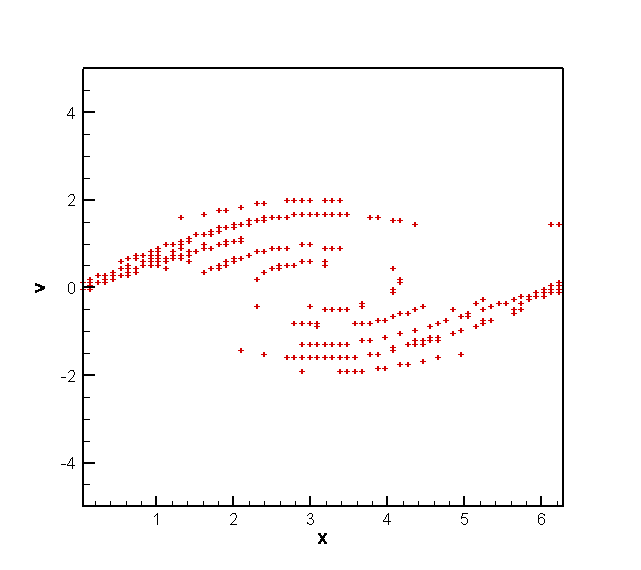} }
\caption{Phase space plots of the two stream instability at $T=53$. The numerical mesh is $64\times128$. Top left: CSLWENO. Top right: CSLHWENO-WO ($CFL=1.2$). Bottom left: CSLHWENO-WL ($CFL=2.2$); the TVB constants $M_x=1,M_y=10$. Bottom right: trouble cells of CSLHWENO-WL at the last time step.}
\label{two}
\end{figure}

\end{example}

\begin{example}
Two stream instability \cite{umeda2008conservative,crouseilles2010conservative}.

Consider the symmetric two stream instability,
\begin{equation*}
f(x,v,t=0) = \frac{1}{2v_{th}\sqrt{2\pi} }\left[ \exp\left(-\frac{(v-u)^2}{2v_{th}^2} \right) + \exp\left(-\frac{v+u}{2v_{th}^2} \right) \right](1+0.05 \cos(kx) )
\end{equation*}
with $u=0.99,\ v_{th}=0.3$ and $k=\frac{2}{13}$.  Our numerical simulation parameters are $v_{max}=5,\ N_x=512,\ N_v=512$.
 In the first row of Figure \ref{two2norm}, the time evolution of the $L^2$ and $L^\infty$ norms of the electric field (in semi-log scale) for CSLWENO and CSLHWENO-WO are plotted.
Figure \ref{two2} shows  numerical solutions of phase space profiles for CSLWENO, CSLHWENO-WO and CSLHWENO-WL at $T=70$. The TVB constants of CSLHWENO-WL are $M_x=M_y=0.1$.

\begin{figure}[h!]
\centering
\includegraphics[height=65mm]{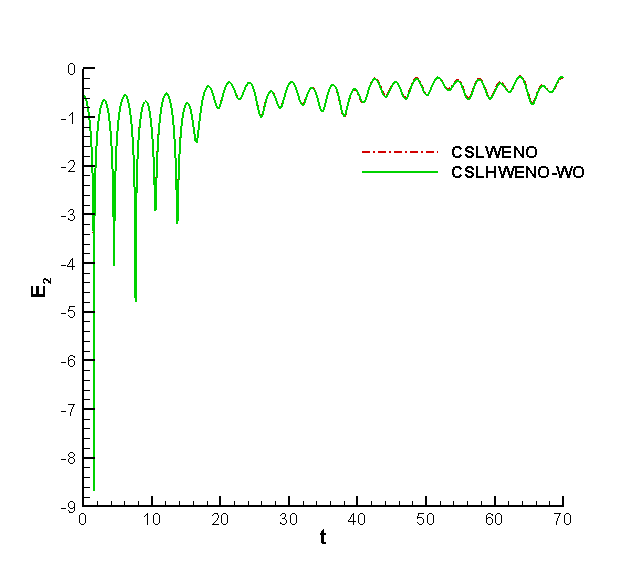}
\includegraphics[height=65mm]{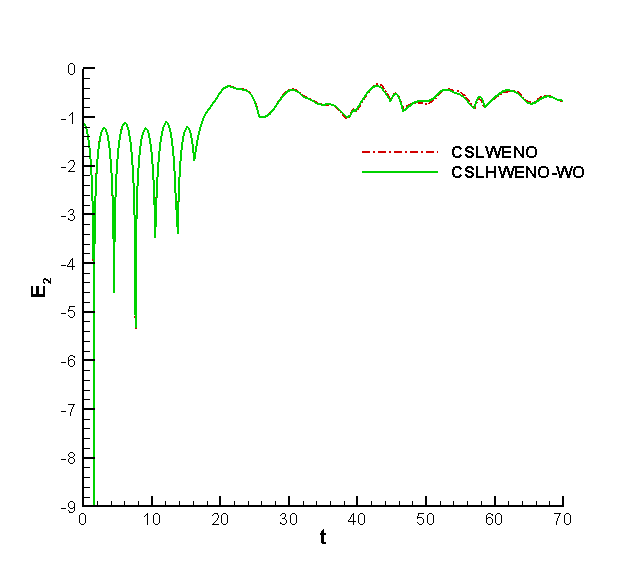}
\includegraphics[height=65mm]{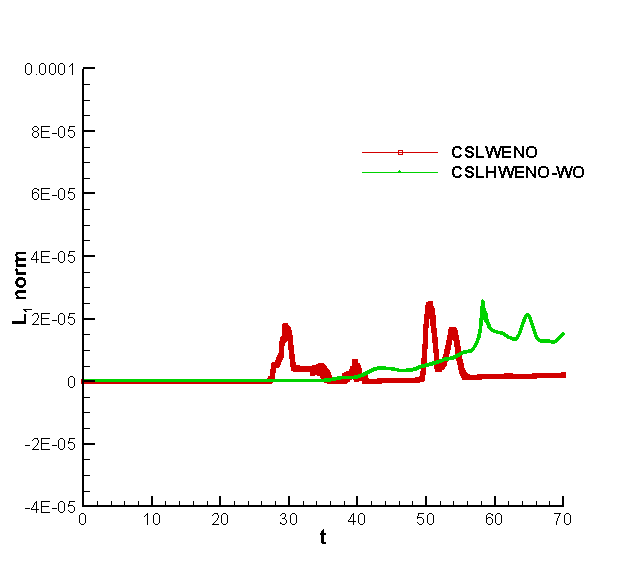}
\includegraphics[height=65mm]{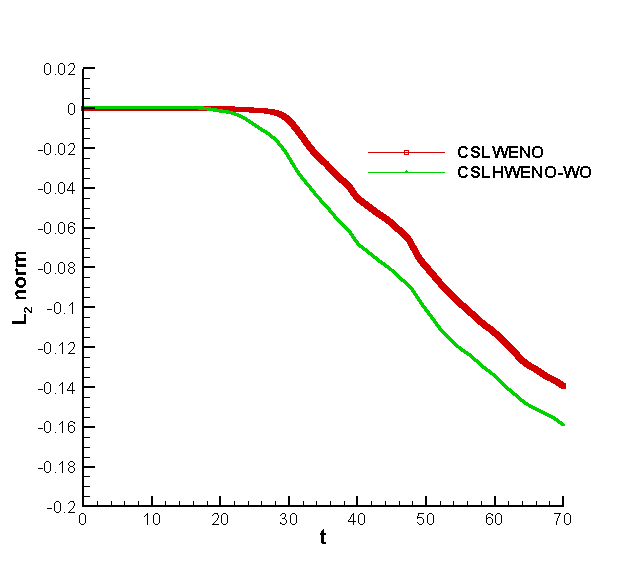}
\includegraphics[height=65mm]{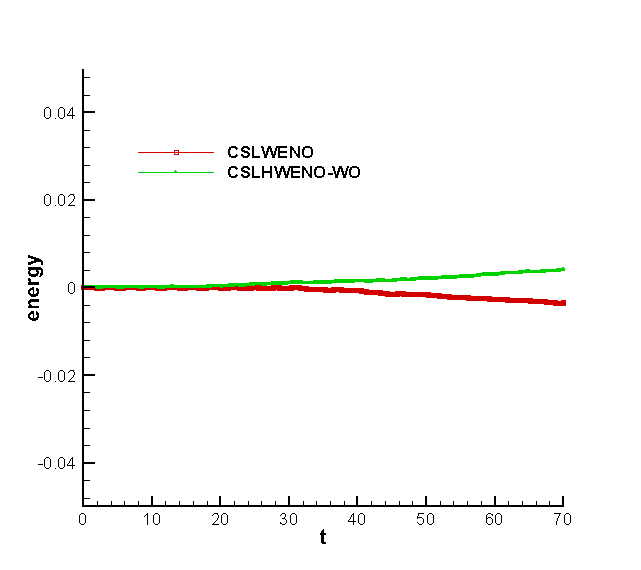}
\includegraphics[height=65mm]{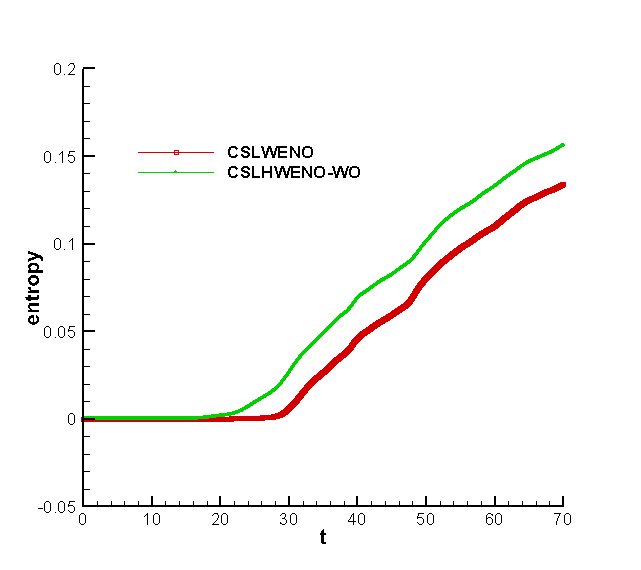}
\caption{Two-stream instability: time evolution of the electric field in $L^2$ (upper left) and $L^\infty$ (upper right) norms, $L^1$ (middle left) and $L^2$ (middle right) norms of the solution as well as the discrete kinetic energy (lower left) and entropy (lower right).}
\label{two2norm}
\end{figure}

\begin{figure}[h!]
\centering
\subfigure[CSLWENO]{
\includegraphics[height=70mm]{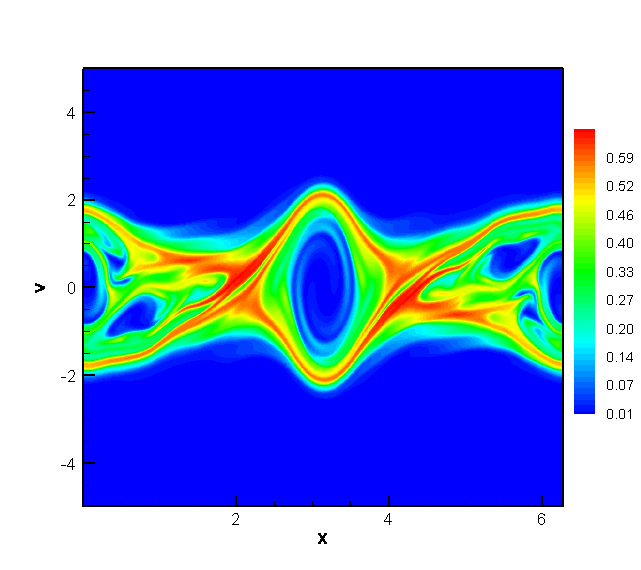} }
\subfigure[CSLHWENO-WO]{
\includegraphics[height=70mm]{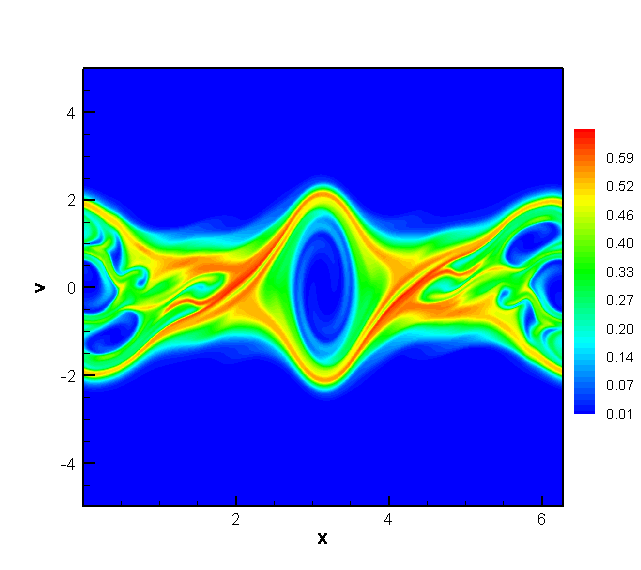} }
\subfigure[CSLHWENO-WL]{
\includegraphics[height=70mm]{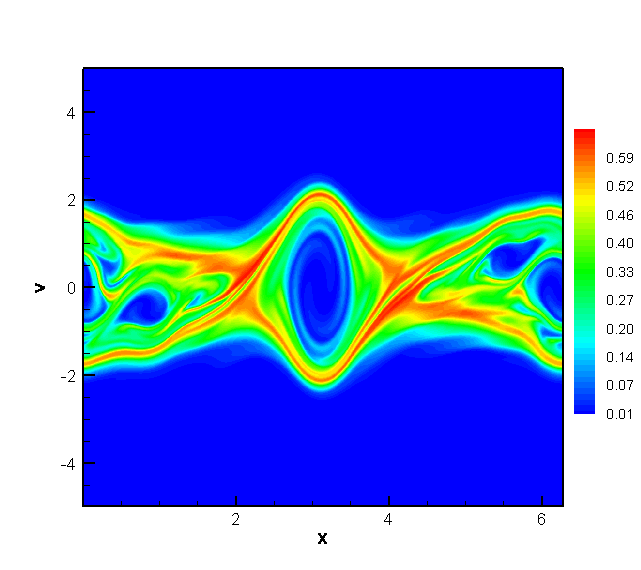} }
\subfigure[trouble cells]{
\includegraphics[height=70mm]{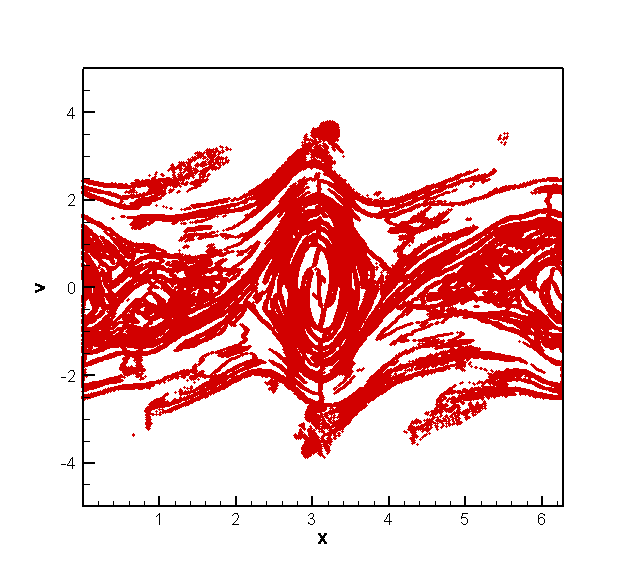} }
\caption{Phase space plots of the two stream instability at $T=70$. The numerical mesh is $512\times512$. Top left: CSLWENO. Top right: CSLHWENO-WO ($CFL=1.2$). Bottom left: CSLHWENO-WL ($CFL=2.2$). Bottom right: trouble cells of CSLHWENO-WL at the last time step.}
\label{two2}
\end{figure}

\end{example}

\section{Conclusions}
\label{conclusion}
In this paper, we propose a conservative SL HWENO scheme for VP system based on dimensional splitting.
Compared with the original WENO reconstruction, the advantage of HWENO reconstruction is compact.
To ensure local mass conservation, the derivative in the scheme is rewritten as the flux-difference form.
The fifth order conservative SL HWENO scheme for the flux difference is proposed.
The scheme can be extended to solve high dimensional problem by the Strang splitting method.
We show the SL HWENO scheme with the Eulerian CFL condition perform well for the classical
Landau damping and the two-steam instability in plasma physics.
When the time stepping size is larger than the Eulerian CFL restriction, we introduce WENO limiters to control oscillations.


\end{document}